\documentclass[10pt]{article}
\usepackage{amssymb,amsbsy,amsmath,amsfonts,amssymb,amscd, mathrsfs}

\usepackage[english]{babel}
\usepackage[T1]{fontenc}
\usepackage{indentfirst}

\makeatletter
\@addtoreset{equation}{section}
\makeatother

\newtheorem{statement}{}[section]
\newtheorem{theoreme}[statement]{Theorem}
\newtheorem{lemme}[statement]{Lemma}

\newtheorem{proposition}[statement]{Proposition}
\newtheorem{definition}[statement]{Definition}

\newcommand\C{\mathbb C}
\newcommand\N{\mathbb N}
\newcommand\R{\mathbb R}
\newcommand\T{\mathbb T}

\newcommand\Z{\mathbb Z}
\newcommand\E{\mathbb E}
\newcommand\e{{\rm e}}
\renewcommand\P{\mathbb P}
\newcommand\eps{\varepsilon}

\let\hat=\widehat
\newcommand\ind{{\rm 1\kern-.30em I}}
\newcommand\qed{\hfill $\square$}
\renewcommand \Re{{\mathfrak R}{\rm e}\,}

\let\phi=\varphi

\newcommand\lpnorm{[\kern-.15em [}
\newcommand\rpnorm{]\kern-.15em ]}
\newcommand\prob{\P}
\newcommand\esp{\E}

\title{\bf Some new thin sets of integers in Harmonic Analysis}
\author{\it Daniel Li,\\ \it Herv\'e Queff\'elec, Luis Rodr{\'\i}guez-Piazza}

\date{}

\begin{document}

\maketitle

{\small\noindent{\bf Abstract}. {\it We randomly construct various subsets $\Lambda$ of the integers which have 
both smallness and largeness properties. They are small since they are very close, in various meanings, to Sidon sets: 
the continuous functions with spectrum in $\Lambda$ have uniformly convergent series, and their Fourier 
coefficients are in $\ell_p$ for all $p>1$; moreover, all the Lebesgue spaces $L^q_\Lambda$ are equal  for 
$q<+\infty$. On the other hand, they are large in the sense that they are dense in the Bohr group and that the 
space of the bounded functions with spectrum in $\Lambda$ is non separable. So these sets are very different 
from the thin sets of integers previously known.\par}}
\medskip

{\small\noindent{\bf R\'esum\'e}. {\it On construit al\'eatoirement des ensembles $\Lambda$ d'entiers positifs 
jouissant simultan\'ement de propri\'et\'es qui les font appara{\^\i}tre \`a la fois comme petits et comme grands. Ils 
sont petits car tr\`es proches \`a plus d'un \'egard des ensembles de Sidon: les fonctions continues \`a spectre dans 
$\Lambda$ ont une s\'erie de Fourier uniform\'ement convergente, et ont des coefficients de Fourier dans $\ell_p$ 
pour tout $p>1$; de plus, tous les espaces de Lebesgue $L^q_\Lambda$ co{\"\i}ncident pour $q<+\infty$. Mais ils 
sont par ailleurs grands au sens o\`u ils sont denses dans le compactifi\'e de Bohr et o\`u l'espace des fonctions 
born\'ees \`a spectre dans $\Lambda$ n'est pas s\'eparable. Ces ensembles sont donc tr\`es diff\'erents des 
ensembles minces d'entiers connus auparavant.\par}}
\medskip

\noindent{\bf Key-words}. ergodic set -- lacunary set -- $\Lambda(q)$-set -- quasi-independent set -- random set -- 
$p$-Rider set -- Rosenthal set -- $p$-Sidon set -- set of uniform convergence -- uniformly distributed set.\par
\medskip

\noindent {\bf Mathematics Subject Classification}. {\it Primary}: 42A36 -- 42A44 -- 42A55 -- 42A61 -- 43A46;   
\noindent {\it Secondary}: 60D05
\vskip 20pt


{\Large\bf Introduction}\medskip

It is well known that the Fourier series of an integrable function defined on the unit-circle $\T=\R/2\pi\Z$ of the 
complex plane $\C$ can be badly behaved. For example, it is well known that there exist continuous functions 
whose Fourier series is not everywhere convergent (see \cite{Korner}, Th. 18.1, and Th. 19.5 for the optimal result), 
and integrable ones with everywhere divergent Fourier series (see \cite{Korner}, Th. 19.2 for instance; see also 
\cite{Konyagin}).\par\smallskip

The problem of thin sets of integers is the following: instead of considering all the integrable functions on 
$\T$, or all the continuous ones, we consider only those whose spectrum (the set where their Fourier coefficients do 
not vanish) is contained in a prescribed subset $\Lambda$ of the integers $\Z$. This set $\Lambda$ will be said 
``thin'' if the Fourier series of these functions behaves better than in the general case. A typical example is 
$\Lambda=\{1,3,3^2,\ldots, 3^n,\ldots\,\}$. It is well known (see \cite{Zyg}, for instance) that every integrable 
function $f$ with spectrum in $\Lambda$ ($f\in L^1_\Lambda$) is actually square integrable, and 
that every continuous function $f$ with spectrum in $\Lambda$ ($f\in{\cal C}_\Lambda$) has a normally 
convergent Fourier series (equivalently $\hat f\in\ell_1$).\par\smallskip

In his seminal paper \cite{Rudin}, W. Rudin defined two notions of thinness for $\Lambda$: $\Lambda$ is a Sidon 
set if $f\in{\cal C}_\Lambda$ implies that $\hat f\in\ell_1$, and $\Lambda$ is a $\Lambda(q)$-set for some 
$q>1$, if $f\in L^1_\Lambda$ implies that $f\in L^q$. These concepts may as well be defined in 
the more general setting of a compact abelian group $G$ equipped with its normalized Haar measure, and for a 
subset $\Lambda$ of its discrete dual group $\Gamma$.\par 
W. Rudin studied the general properties of those sets and the connection between the two notions. In particular, he 
showed that Sidon sets are $\Lambda(q)$-sets for all $q<+\infty$, and that, more precisely:

\begin{itemize}
\item [(0.1)] $\Lambda$ Sidon implies $\Vert f \Vert_q\leq C\sqrt q \Vert f \Vert_2$ for every 
$\Lambda$-polynomial $f$ and for every $q\geq 2$, where $C$ is a constant which depends only on the Sidon 
constant of $\Lambda$.
\end{itemize}

Since then, several new notions of thin sets emerged. These include $p$-Sidon sets (see 
\cite{Blei1}, \cite{Blei2}, \cite{Blei3}, \cite{Blei-K}, \cite{Bo-Pyt}, \cite{Ed-Ross}, \cite{Fournier-P}, 
\cite{Johnson}, \cite{Johnson-W}, \cite{Lefevre2}, \cite{Pisier2}, \cite{Woodward}), and sets of uniform 
convergence (see \cite{Arki-Oskol}, \cite{Figa}, \cite{Fournier}, \cite{Fournier-P}, \cite{Kashin-T}, 
\cite{Lefevre1}, \cite{Oskol}, \cite{Pedem}, \cite{Soardi-T}, \cite{Trava}): every continuous function with spectrum 
in such a set has its Fourier seriesin $\ell_p$ or uniformly convergent,  respectively. But the examples of such 
sets were always nearly the same: products (sometimes ``fractional products'': \cite{Blei2}, \cite{Blei3}, 
\cite{Blei-K}), or sums of Sidon sets, which is a severe restriction for the geometry of the Banach space 
${\cal C}_\Lambda$. For example, F. Lust--Piquard (\cite{Lust1}) proved that:

\begin{itemize}
\item [(0.2)] The injective tensor product $\ell_1\hat\otimes_\eps\cdots\hat\otimes_\eps\ell_1$ has the Schur 
property ({\it i.e.} weakly null sequences converge in norm to zero).
\end{itemize}

It follows easily that:

\begin{itemize}
\item [(0.3)] If $\Lambda=E_1\times\cdots\times E_k$, where the $E_j$'s are Sidon sets, then 
${\cal C}_\Lambda$ has the Schur property; in particular, ${\cal C}_\Lambda$ does not contain $c_0$, the space 
of sequences going to zero at infinity.
\end{itemize}

Since these sets were essentially the only known examples of $p$-Sidon sets (they are exactly $2N/(N+1)$-Sidon), 
one could believe that all $p$-Sidon sets have this property. It should be mentionned that in \cite{Blei2}, R. Blei 
constructed for each $p\in]1,2[$, exactly $p$-Sidon sets, using fractional products, so of a different type, but the 
corresponding space ${\cal C}_\Lambda$ appears as an $\ell_1$-sum of finite dimensional spaces, and so does 
have the Schur property (we thank R. Blei for this remark).\par

Because of this lack of examples, the comparison between two classes of thin sets proved to be very difficult: whether 
a $p$-Sidon, or a set of uniform convergence is a $\Lambda(q)$-set for some $q>1$ is still an open problem. On the 
other hand, considerable progress concerning the Sidon sets or $\Lambda(q)$-sets has been made: for example, 
G. Pisier (\cite{Pisier2}, Th. 6.2) proved that the converse of (0.1) is true, and J. Bourgain (\cite{Bourgain12}) 
proved that for each $q>2$ there exist ``exactly'' $\Lambda(q)$-sets, {\it i.e.} sets which are $\Lambda(q)$, but 
$\Lambda(q')$ for no $q'>q$. Both authors used random methods, and more specifically, J. Bourgain popularized 
the ``method of selectors'' to produce several thin sets $\Lambda$ with unusual properties, such as being 
``uniformly distributed'', which implies, by a result of F. Lust--Piquard (\cite{Lust42}), that ${\cal C}_\Lambda$ 
contains $c_0$ and therefore is not a Rosenthal set ({\it i.e.} there are bounded measurable functions with 
spectrum in $\Lambda$ wich are not almost everywhere equal to a continuous function), and which also implies 
that $\Lambda$ is dense in the Bohr group (see \cite{Blum}, Theorem 1). This allowed the first named author to 
see that there are sets of integers which are $\Lambda(q)$ for all $q<+\infty$ but not Rosenthal (\cite{Li}; see also 
\cite{Neuwirth}).\par
\medskip

The aim of this paper is the construction of random sets $\Lambda$ of integers which have thinness properties, but 
wich are not Rosenthal sets ({\it i.e.} ${\cal C}_\Lambda$ is not the whole $L^\infty_\Lambda$), actually such 
that ${\cal C}_\Lambda$ contains $c_0$, and are dense in the Bohr group. In view of (0.3), these sets will 
necessarily be very exotic compared to the previously known examples. This shows that replacing absolute 
convergence of the Fourier series by uniform convergence (sets of uniform convergence) or by $\ell_p$ convergence 
for $p>1$ ($p$-Sidon sets) gives sets which are very far from Sidon sets. This constrasts 
with Pisier's result saying that $\Lambda$ is necessarily a Sidon set whenever $\hat f\in \ell_{1,\infty}$ for every 
$f\in {\cal C}_\Lambda$ (from \cite{Pisier48}, Th\'eor\`eme 2.3 (vi), and the top of page 688). On the other hand, 
though non-Sidon Rosenthal sets do exist (\cite{Rosenthal}), it follows from Bourgain-Milman's cotype theorem 
(\cite{Bourgain-Milman}) that, for every non-Sidon set $\Lambda$, ${\cal C}_\Lambda$ does contain 
$\ell_\infty^n$ uniformly, so that the presence of $c_0$ inside ${\cal C}_\Lambda$ for non-Sidon $\Lambda$ 
may appear not so surprising. Although it is not known whether Sidon sets may be dense in 
the Bohr group, we obtain in this paper, as mentioned above, sets which are dense in the Bohr group, and are of 
uniform convergence and $p$-Sidon for every $p>1$.\par
\smallskip

We construct essentially four types of sets. Each of them will be a non Rosenthal set, but a set of uniform 
convergence, $\Lambda(q)$ for all $q<+\infty$, and with moreover additional properties of $p$-Sidonicity.\par 
The first one (Theorem~\ref{theo 2.2}) is a very lacunary set $\Lambda$ with the nicest properties: 
it is $p$-Sidon for all $p>1$. The second and third ones (Theorem~\ref{theo 2.5} and Theorem~\ref{theo 2.6}) 
are medium lacunary sets: for each $p$ with $1<p<4/3$, they are, in Theorem~\ref{theo 2.5}, $p$-Rider 
(a weaker property than being $p$-Sidon, see the definition below), but not $q$-Rider for $q<p$, and are $q$-Sidon 
for every $q>p/(2-p)\,;$ and in Theorem~\ref{theo 2.6}, they are $q$-Rider for every $q>p$, but not $p$-Rider, 
and they are $q$-Sidon for every $q>p/(2-p)$. Finally, the fourth type (Theorem~\ref{theo 2.7}) 
is a set $\Lambda$ which is, in some sense as little lacunary as possible if we want its trace on each interval 
$[N, 2N[$ to have a bounded Sidon constant. It leads to sets which are  4/3-Rider, but not  $q$-Rider for $q<4/3$.
\par

We construct these sets by using various choices of selectors, and adding arithmetical, functional or probabilistic 
arguments. The treatment of the last case requires a different probabilistic approach, taken from \cite{Bourgain8}.
\par
It should be noted that in the two first cases the sets are uniformly distributed; in the fourth case , however, the sets 
$\Lambda$ only have positive upper density in uniformly distributed sets. Nevertheless, ${\cal C}_\Lambda$ still 
contains $c_0$, by a result of F. Lust-Piquard (\cite{Lust42}, Th. 5).\par
\bigskip

\noindent{\bf Acknowledgement.} Part of this paper was made when the first named author was a guest of the 
Departamento de An\'alisis Matem\'atico de la Universidad de Sevilla in April 1999, and when the third named 
author was a guest of the Universit\'e d'Artois in Lens in june 1999.


\section{Notation, definitions and preliminary results}

We denote by $\T$ the compact abelian group of complex numbers of modulus one, equipped with its normalized 
Haar measure $m$. ${\cal C}(\T)$ denotes the space of continuous complex functions defined on $\T$, equipped 
with its $\sup$ norm $\Vert\ \Vert_\infty$ and identified as usual with the space of continuous $2\pi$-periodic 
complex functions defined on $\R$. If $\Lambda$ is a subset of the dual group $\Z$, ${\cal C}_\Lambda$ will 
denote the subspace of ${\cal C}(\T)$ consisting of functions whose spectrum lies in $\Lambda$: 
$$\hat f(n)\equiv\int_\T fe_{-n}\,dm=0
\hbox{\hskip 5mm if } n\in\Z\setminus\Lambda,$$
where $e_n(z)=z^n$, or equivalently, $e_n(t)=\e^{int}$.\par

${\cal C}_\Lambda$ is the uniform closure of the space ${\cal P}_\Lambda$ of trigonometric polynomials with 
spectrum in $\Lambda$, {\it i.e.} the uniform closure of the subspace ${\cal P}_\Lambda$ generated by the 
characters $e_n$, with $n\in\Lambda$.\par

For $f\in {\cal C}(\T)$, $1\leq q<+\infty$, $M$ and $N$ positive integers, we shall denote the Fourier sums of 
$f$ by:
$$S_{M,N}(f)=\sum_{-M}^N\hat f(n)e_n$$
and the symmetric Fourier sums of $f$ by:
$$S_N(f)=S_{N,N}(f)=\sum_{-N}^N\hat f(n)e_n\,.$$\par

$\vert A \vert$ denotes the cardinality of the finite set $A$.\par
\smallskip

A {\it relation} in $\Lambda\subseteq\Z^\ast\equiv\Z\setminus\{0\}$ is a $(+1,-1,\,0)$-valued sequence 
$(\theta_k)_{k\in\Lambda}$ such that $\sum\vert\theta_k\vert<+\infty$ and $\sum\theta_k k=0$. The set 
$S=\{k\,;\ \theta_k\not=0\}$ is called the {\it support} of the the relation, and 
$\vert S \vert=\sum\vert \theta_k\vert$ is called its {\it length}.\par

The relation $(\theta'_k)_{k\in\Lambda}$ is said to be {\it longer} than the relation $(\theta_k)_{k\in\Lambda}$ 
if $\theta_k\not=0$ implies $\theta_k=\theta'_k$.\par

The set $\Lambda\subseteq \Z^\ast$ is {\it quasi-independent} if it contains no non-trivial relation ({\it i.e.} with 
non-empty support). Typically, $\Lambda=\{1,2,4,\ldots,2^n,\ldots\,\}$ is quasi-independent. The 
quasi-independent sets are the prototype of Sidon sets, {\it i.e.} of sets $\Lambda$ for which: 
$\Vert \hat f \Vert_1\leq K\Vert f \Vert_\infty$ for all $f\in{\cal C}_\Lambda$. The best constant $K$ in this 
inequality is called the Sidon constant of $\Lambda$ and is denoted by $S(\Lambda)$. We will refer to 
\cite{Lopez-Ross} for standard notions on Sidon sets. It is known that quasi-independent sets are not only Sidon sets 
but their Sidon constant is bounded by an absolute constant: this follows from \cite{Rudin}, Th. 2.4 and 
\cite{Rider1}, Lemma 1.7. Other proofs can be found in \cite{Pisier48}, lemme 3.2, and in 
\cite{Bourgain9}, Prop. 1. We shall use the fact that $S(\Lambda)\leq 8$ if $\Lambda$ is quasi-independent.\par
\smallskip

Let us recall now some classical definitions and results.\par
\smallskip

A set $\Lambda\subseteq\Z$ is said to be a {\it $\Lambda(q)$-set} (where $q>2$) if there exists a positive 
constant $C_q$ such that $\Vert f\Vert_q\leq C_q\Vert f\Vert_2$ for every $f\in{\cal P}_\Lambda$.\par
The notion of a $\Lambda(q)$-set is, in some sense, local. That follows from the Littlewood-Paley theory. The next 
proposition is essentially well-known, except for the growth of the constant, for which we have found no reference. 
Accordingly, we offer a short proof.\par

\begin{proposition}\label{prop 1.1}
Let $\Lambda\subseteq[2,+\infty[$. Then:\par
{\rm (a)} Let $(M_n)_{n\geq1}$ be a sequence of positive integers such that $M_1\leq 2$ and 
$M_{n+1}/M_n\geq\alpha>1$. If $\Lambda\cap[M_n,M_{n+1}[$, $n\geq1$, has a uniformly bounded Sidon 
constant, then $\Lambda$ is $\Lambda(q)$ for all $q\geq2$; more precisely: 
$\Vert f\Vert_q\leq C(q,\alpha)\Vert f\Vert_2$ for every $f\in{\cal P}_\Lambda$.\par
{\rm (b)} If $\Lambda\cap[2^n,2^{n+1}[$, $n\geq1$, has a uniformly bounded Sidon constant, $\Lambda$ is 
$\Lambda(q)$ for every $q\geq2$ and, more precisely: $\Vert f\Vert_q\leq Cq^2\Vert f\Vert_2$ for every 
$f\in{\cal P}_\Lambda$ and for some numerical constant $C$.
\end{proposition}

\noindent{\bf Proof.} (a) Set 
$$f_k=\sum_{M_k\leq n<M_{k+1}} \hat f(n)e_n 
\qquad \text{and} \qquad Sf=\Big(\sum_{k=1}^{+\infty}\vert f_k\vert^2\Big)^{1/2}.$$ 
Since $M_{k+1}/M_k\geq\alpha>1$ and $\Lambda\subseteq[M_1,+\infty[$, we have 
(\cite{Zyg}, Chap. XV, Th. 2.1):
$$\Vert f\Vert_q\leq C_0(q,\alpha)\Vert Sf\Vert_q\,.$$
Now, using the 2-convexity of the $L^q$-norm for $q\geq2$, we obtain:
$$\Vert Sf\Vert_q\leq\Big(\sum\limits_{k=1}^{+\infty}
\Vert f_k\Vert_q^2\Big)^{1/2}.$$ 
But $f_k\in{\cal P}_{\Lambda_k}$, where $\Lambda_k=\Lambda\cap[M_k,M_{k+1}[$ has a uniformly bounded 
Sidon constant. Therefore $\Vert f_k\Vert_q\leq C_1\sqrt q \Vert f\Vert_2$, where $C_1$ is a numerical 
constant. The result follows.\par
(b) We now make use of the classical square function 
$$Sg=\Big(\sum\limits_{k\in\Z} \vert g_k\vert^2\Big)^{1/2},$$ 
where
$$g_k=\sum_{2^k\leq n<2^{k+1}}\hat g(n)e_n\,,\ {\rm if}\ k\geq0 \hskip 5mm 
{\rm and}\hskip 5mm 
g_k=\sum_{-2^{\vert k\vert+1}< n\leq -2^{\vert k\vert}}\hat g(n)e_n\,
\ {\rm if}\ k<0\ .$$\par

For this classical square function, we have the following sharp inequality, due to J. Bourgain 
(\cite{Bourgain11}, Th. 1):
$$\Vert Sg\Vert_p\leq C_0 (p-1)^{-3/2}\Vert g\Vert_p\hskip 5mm{\rm for}\ 
1<p\leq2\,,$$
where $C_0$ is a numerical constant. We deduce by duality that:
$$\Vert f\Vert_q\leq C_0 q^{3/2}\Vert Sf\Vert_q\hskip 5mm{\rm for}\ 
2\leq q<+\infty\,.$$
In fact, by orthogonality (recall that $f\in{\cal P}_\Lambda$ and that $\Lambda\subseteq[2,+\infty[$) and 
the Cauchy--Schwarz inequality, we have, for every $g\in L^p$ with $\Vert g\Vert_p=1$ ($1/p+1/q=1$):
\begin{align*}
\vert<f,g>\vert=&\big\vert\sum_{k=1}^{+\infty}<f_k,g_k>\big\vert=
\big\vert\int_\T\ \sum_{k=1}^{+\infty} f_k(-t)g_k(t)\,dm(t)\big\vert \\
& \leq\int_\T Sf(-t)Sg(t)\,dm(t)\cr
&\leq\Vert Sf\Vert_q\Vert Sg\Vert_p\leq 
C_0 (p-1)^{-3/2}\Vert Sf\Vert_q \\
& \leq C_0 q^{3/2}\Vert Sf\Vert_q.
\end{align*}
\par

This means that here we are allowed to take $C_0(q,2)=C_0q^{3/2}$ in part (a) of the proof. The rest is unchanged, 
and we can also take $C(q,2)=C_1\sqrt q C_0\, q^{3/2}= C q^2$.\hfill\qed
\bigskip

A set $\Lambda\subseteq\Z$ is called a {\it set of uniform convergence} (in short a {\it UC-set}) if, for any 
$f\in{\cal C}_\Lambda$, the symmetric Fourier sums $S_N(f)$ converge uniformly to $f$. Its constant of uniform 
convergence $U(\Lambda)$ is the smallest constant $K$ such that, for any $f\in{\cal C}_\Lambda$:
$$\sup_N\Vert S_N(f)\Vert_\infty\leq K\Vert f\Vert_\infty.$$\par

The following variant turns out to be more tractable ([56]). $\Lambda$ is called a {\it set of complete uniform 
convergence} (in short a CUC-set) if the translates ($\Lambda+a$) are uniformly UC for $a\in\Z$, or equivalently, 
if the Fourier sums $S_{M,N}(f)$ converge uniformly to $f$ as $M,N$ go to $+\infty$, for every 
$f\in {\cal C}_\Lambda$.\par
The two notions turn out to be distinct (\cite{Fournier}), but clearly coincide if $\Lambda\subseteq\N$, which 
will always be the case in the sequel. The notion of CUC-set is also a local one as the following proposition shows.

\begin{proposition}[\cite{Trava}, Th. 3]\label{prop 1.2}
Let $\Lambda\subseteq\N^\ast$ and $\Lambda_N=\Lambda\cap[N,2N[$.\par
{\rm (a)} If $U(\Lambda_N)$ is bounded by $K$ for $N=1,2,\ldots$, then $\Lambda$ is a {\rm CUC}-set.\par
{\rm (b)} Let $(M_n)_{n\geq1}$ be a sequence of positive integers such that $M_{n+1}/ M_n\geq2$. Then, if
$\Lambda\cap[M_n,M_{n+1}[$ are quasi-independent for each $n$, or more generally if they are Sidon sets with 
uniformly bounded Sidon constant, then $\Lambda$ is a {\rm CUC}-set.
\end{proposition}

\noindent{\bf Remark.} (b) is a useful criterion to produce sets that are CUC but not Sidon; for instance, if 
$\Lambda=\bigcup_{n=1}^{+\infty}\{2^n+2^j\,;\ j=0,\ldots,n-1\,\}$, then $\Lambda\cap[2^n,2^{n+1}[$ is 
quasi-independent, whereas $\Lambda\cap[1,N]$ has about $(\log N)^2$ elements, and therefore cannot be 
Sidon (the mesh condition for Sidon sets, see Proposition~\ref{prop 1.6} below, is violated).\par
\bigskip

The random variables which we shall use will always be defined on some probability space 
$(\Omega,{\cal A},\P)$ which will play no explicit role, and the expectation with respect to $\P$ will always be 
denoted by $\esp$:
$$\esp(X)=\int_\Omega X(\omega)\,d\prob(\omega)\ .$$\par

Recall the (more or less) classical deviation inequality (see \cite{Ledoux-Ta}, \S~6.3):

\begin{lemme}\label{lemme 1.3}
Let $X_1,\ldots,X_N$ be independent centered complex random variables such that $\vert X_k\vert\leq 1$, 
$k=1,\ldots,N$. Let $ \sigma\geq\sum\limits_{k=1}^N \esp\vert X_k\vert^2$. Then, one has, for every 
$a\leq\sigma$:
$$\prob(\vert X_1+\cdots+X_N\vert\geq a)\leq  4\exp (- a^2 / 8\sigma)\,.$$
\end{lemme}

Let $(r_n)_n$ be a Bernoulli sequence, {\it i.e.} a sequence of independent random variables such that:
$$\prob(r_n=1)=\prob(r_n=-1)= 1/2\, .$$\par

For $f\in {\cal P}$, the space of trigonometric polynomials, $\lpnorm f\rpnorm$ denotes the norm of $f$ in the 
Pisier's space ${\cal C}^{\hbox{\it a.s.}}$:
$$\lpnorm f\rpnorm=\esp\big\Vert \sum_n r_n \hat f(n)e_n\big\Vert_\infty\,.$$

See \cite{Kahane} and \cite{Pisier2} for more information about this norm.\par

\begin{definition}\label{def 1.4}
A set $\Lambda\subseteq\Z$ is called a {\it $p$-Sidon set} ($1\leq p<2$) if there exists a constant $K$ such that 
$\Vert \hat f \Vert_p\leq K\Vert f \Vert_\infty$ for all $f\in {\cal P}_\Lambda$.\par
It is said to be a {\it $p$-Rider set} if there exists a constant $K$ such that 
$\Vert \hat f \Vert_p\leq K\lpnorm f\rpnorm$ for all $f\in {\cal P}_\Lambda$.
\end{definition}

$p$-Rider sets were implicitely introduced, with different definition, in \cite{Ed-Ross} (Th. 2.4), and in 
\cite{Johnson}, p. 213, as class ${\cal T}_p$ (see also \cite{Pisier2}, Th. 6.3). They were explicitely defined 
and studied in \cite{Luis1} and \cite{Luis2} under the name ``$p$-Sidon presque s\^urs''. We used 
``almost surely $p$-Sidon set'' in the first version of this paper, but, following a suggestion of J.-P. Kahane, we 
now use the terminology ``$p$-Rider''.\par
Clearly, every $p$-Sidon set is $p$-Rider. The converse is true for $p=1$: this is a remarkable result due to 
D. Rider (\cite{Rider2}), making clever use of Drury's convolution device (which proves that the union 
of two Sidon sets is Sidon \cite{Drury}). Whether this converse is still true for $1<p<2$ is an open problem.\par

\begin{definition}\label{def 1.5} 
We shall say that a finite set $B\subseteq \Lambda$ is $M$-{\it pseudo-complemented} in $\Lambda$ if there 
exists a measure $\mu$ on $\T$ such that:
$$\vert\hat\mu\vert\geq1\hbox{ on }B\,;\hskip 3mm 
\hat\mu=0\hbox{ on }\Lambda\setminus B\,;\hskip 3mm 
\Vert\mu\Vert\leq M\,.$$
\end{definition}

The following proposition gives some necessary, sufficient, or necessary and sufficient conditions for a set 
$\Lambda$ to be $p$-Sidon or  $p$-Rider. Part (b) of this proposition seems to be new.

\begin{proposition}\label{prop 1.6}
Let $\Lambda\subseteq\Z^\ast$ and $1\leq p<2$. Set $\eps(p)=2/p-1$. Then:\par
{\rm (a)} $\Lambda$ is a $p$-Rider set if and only if there exists a constant $\delta>0$ such that, for every finite 
set $A\subseteq\Lambda$, there exists a quasi-independent subset $B\subseteq A$ such that 
$\vert B\vert\geq \delta\vert A\vert^{\eps(p)}$.\par
{\rm (b)} Let $q_0>1$. If there exists a constant $\delta>0$ such that, for every finite set $A\subseteq\Lambda$, 
there exists a quasi-independent subset $B\subseteq A$ such that $\vert B\vert\geq \delta\vert A\vert^{1/q_0}$ 
and if $B$ can moreover be taken $M$-pseudo-complemented in $\Lambda$, for some fixed $M$, then 
$\Lambda$ is a $q$-Sidon set for every $q>q_0$.\par
{\rm (c)} If $\Lambda$ is a $p$-Rider set, we have the following mesh condition:
$$\vert\Lambda\cap[1,N]\vert\leq C(\log N)^{p/(2-p)}.$$
\end{proposition}

\noindent{\bf Proof.} We refer to \cite{Luis1} for the proof of (a) and (c). To prove (b), let $f\in {\cal P}_\Lambda$, 
fix $t>0$, and set $A=\{\,\vert \hat f\vert>t\,\}$. Take $B\subseteq A$ and $\mu$ as in Definition~\ref{def 1.4}. 
Then $B$ is a Sidon set with Sidon constant $\leq 8$, and since 
$f\ast\mu=\sum\limits_{n\in B} \hat f(n)\hat \mu(n)e_n$,
\begin{align*}
\Vert f \Vert_\infty
\geq M^{-1}\Big\Vert\sum_B \hat f(n) \hat\mu(n) e_n \Big\Vert_\infty 
&\geq \frac{1}{8M} \sum_B \vert \hat f(n)\vert\,\vert\hat \mu (n)\vert \cr
& \geq \frac{1}{8M}\sum_B \vert\hat f(n)\vert\geq \frac{t\vert B\vert}{8M}
\geq \frac{t\delta\vert A\vert^{1/q_0}} {8M}\cdot
\end{align*}

In other words, for some constant $C>0$, one has:
$$t\,.\vert\,\{\,\vert\hat f\vert>t\,\}\,\vert^{1/q_0} \leq C\Vert f\Vert_\infty\,,
\hskip 5mm\hbox{for every }t>0\,,$$
which means that the Lorentz norm of $\hat f$ in the Lorentz space $\ell_{q_0,\infty}$ is dominated by 
$\Vert f\Vert_\infty$.\par 
Now, $\ell_{q_0,\infty}$ is continuously injected in $\ell_q$ for $q>q_0$ (see for instance 
\cite{L-T}, II p. 143), and this gives the desired result.\hfill\qed
\bigskip

We denote, as usual, by $c_0$ the classical space of sequences $x=(x_n)_{n\geq0}$ tending to zero at infinity, 
equipped with the norm $\Vert x \Vert=\sup_n\vert x_n \vert$. We say, in the usual familiar way, that a 
Banach space $X$ ``contains $c_0$'' if $X$ has a closed subspace isomorphic to $c_0$. Our notation for Banach 
spaces is classical, as can be found in \cite{Diestel}, \cite{L-T} or \cite{Wojt} for instance.\par

A subset $\Lambda$ of $\Z$ is said to be a {\it Rosenthal set} if every bounded measurable function on $\T$ with 
spectrum in $\Lambda$ is almost everywhere equal to a continuous function (in short 
$L^\infty_\Lambda={\cal C}_\Lambda$). $\Lambda$ is not Rosenthal if and only if $L^\infty_\Lambda$ is not 
separable, so such a set can be thought as being a big set.\par
Every Sidon set is clearly Rosenthal, but H.P. Rosenthal gave examples of 
non-Sidon sets which are Rosenthal (\cite{Rosenthal}). We shall make use of the following well known negative 
criterion (see \cite{Lust41}, \S\,3), which follows from the classical theorem of C. Bessaga and A. Pe\l czy\'nski 
(\cite{L-T}, I.2.e.8), saying that a dual space which contains $c_0$ has to contain also $\ell_\infty$.\par

\begin{proposition}\label{prop 1.7}
If ${\cal C}_\Lambda$ contains $c_0$, then $\Lambda$ is not a Rosenthal set.
\end{proposition}

\begin{definition}\label{def 1.8}
Let $\Lambda\subseteq\N^\ast\equiv\N\setminus\{0\}$, and set 
$$\Lambda_N=\Lambda\cap[1,N] \qquad \text{and} \qquad  
A_N(t)={\frac{1}{\vert \Lambda_N\vert}}\sum\limits_{n\in\Lambda_N} e_n(t)\,.$$ 
We say that $\Lambda$ is:\par
- \emph{ ergodic} if $(A_N(t))_{N\geq1}$ converges to a limit $l_\Lambda(t)\in\C$ for each $t\in\T$.\par
- \emph{strongly ergodic} if it is ergodic and moreover the limit function $l_\Lambda$ defines an element of 
$c_0(\T)$: for every $\eps>0$ the set $\{t\in\T\,;\ \vert l_\Lambda (t)\vert>\eps\,\}$ is finite.
\par
- \emph{uniformly distributed} if it is (strongly) ergodic and, moreover, $l_\Lambda (t)=0$ for 
$t\not=0$ mod. $2\pi$.
\end{definition}

The reason for this terminology is that the ergodic sets are those for which an ergodic theorem holds: 
$(1/\vert \Lambda_N\vert) \sum_{n\in \Lambda_N}T^n$ converges in the strong operator topology for 
every contraction $T$ of a Hilbert space. Typically, the set of $d^{th}$ perfect powers, or the set of prime numbers 
are strongly ergodic (according to the result of Vinogradov for $t$ irrational mod. $2\pi$, and to the 
Dirichlet's arithmetic progression theorem for $t$ rational mod. $2\pi$). The third name comes from 
H. Weyl's classical criterion for the equidistribution of a real sequence mod. $2\pi$.\par

The relationship between these notions comes from:

\begin{theoreme}[F. Lust-Piquard \cite{Lust42}]\label{theo 1.9}
Let $\Lambda\subseteq [1,+\infty[$ be a set of positive integers. Then:\par
{\rm (a)} If $\Lambda$ is strongly ergodic, ${\cal C}_\Lambda$ contains $c_0.$\par
{\rm (b)} More generally, if $\Lambda$ is strongly ergodic and $D\subseteq\Lambda$ has a positive upper density 
with respect to $\Lambda$, then ${\cal C}_D$ contains $c_0$ as well.
\end{theoreme}

Here ``positive upper density'' means that:
$$\mathop{\overline{\rm lim}}_{N\to+\infty} \frac{\vert D\cap[1,N]\vert}{\vert \Lambda\cap[1,N]\vert}>0\,.$$
\par

Part (b) will be useful to us in the last theorem of Section~\ref{section 2}.\par
See \cite{Lust42} for the proof of this theorem. The underlying idea for (a) is that if 
$A_N(t)\to l_\Lambda(t)$ for every $t\in\T$, $l_\Lambda$ defines an element of the biorthogonal 
${\cal C}_\Lambda^{\perp\perp}$, and the condition $l_\Lambda\in c_0(\T)$ implies that it is the sum of a 
weakly unconditionally Cauchy series of continuous functions. By using a perturbation argument due to 
A. Pe{\l}czy\'nski (see \cite{Singer}, lemma 15.7, p. 446) and the classical Bessaga-Pe{\l}czy\'nski theorem, one 
obtains that ${\cal C}_\Lambda$ contains $c_0.$\par\smallskip

This theorem allowed its author to prove that ${\cal C}_\Lambda$ contains $c_0$ when 
$\Lambda=\{1,2^d,3^d,\ldots\}$ is the set of the $d^{th}$ perfect powers, and when 
$\Lambda=\{2,3,5,7,\ldots\}$ is the set of the prime numbers. On the other hand, K. I. Oskolkov 
(\cite{Oskol}; see also \cite{Arki-Oskol}) showed that the set of the $d^{th}$ powers is not a UC-set, and  
J. Fournier and L. Pigno (\cite{Fournier-P}, Th. 4) proved that the set of prime numbers is not a UC-set either. This 
could be taken as an indication that containing $c_0$ is an obstruction to being UC. As we shall see in the next 
section, this is far from being the case: there do exist sets $\Lambda$ which are UC and for which 
${\cal C}_\Lambda$ contains $c_0.$\par\medskip

The last ingredient we require is a random procedure to produce ergodic sets.\par
Let $(\eps_k)_{k\geq1}$ be a sequence of independent $0 - 1$ valued random variables, called ``{\it selectors}'' 
according to the terminology coined by J. Bourgain. To those selectors is associated a random set $\Lambda$ of 
positive integers 
$$\Lambda=\Lambda(\omega)=\{k\geq1\,;\ \eps_k(\omega)=1\,\}\,.$$ 

\begin{theoreme}[J. Bourgain \cite{Bourgain10}, Prop. 8.2]\label{theo 1.10}
\hskip -4pt 
Let $\eps_1,\ldots,\eps_N,\ldots$ be selectors of respective expectations $\delta_1,\ldots,\delta_N,\ldots$ and 
assume that $\sigma_N/\log N\to+\infty$, where $\sigma_N=\delta_1+\cdots+\delta_N$ (which is in particular 
the case when $k\delta_k\to+\infty$), and that $(\delta_n)_n$ decreases. Then the set 
$\Lambda=\Lambda(\omega)$ is almost surely uniformly distributed. In particular, it is almost surely strongly 
ergodic.
\end{theoreme}


\section{Main results}\label{section 2}

In this section, we will always consider selectors $\eps_n$, $n\geq1$, with mean $\delta_n=\alpha_n/n$, with 
$(\alpha_n)_n$ tending to infinity and $(\delta_n)_n$ decreasing.\par
Moreover, except in the last theorem of this section, we will assume that $(\alpha_n)_n$ is increasing.\par

If $\Lambda=\Lambda(\omega)=\{n\geq1\,;\ \eps_n(\omega)=1\}$  is the corresponding random set of integers, 
$\Lambda$ is almost surely uniformly distributed by Bourgain's theorem. Moreover, it also has the nice almost 
sure property of being asymptotically independent; more precisely, there exists an increasing sequence 
$(M_n)_n$ of positive integers such that $\Lambda\cap [M_n,+\infty[$ is both large and without relations of 
length $\leq n$. A subset $B$ of $\Lambda\cap [M_n,+\infty[$ with $n$ elements is then automatically 
quasi-independent, and this allows us to use Propositions~\ref{prop 1.1}, \ref{prop 1.2}, \ref{prop 1.6} to show 
that $\Lambda$ has good additional properties: UC, $p$-Sidon, {\it etc...}\,. To obtain this asymptotic  
quasi-independence, the following half-combinatorial, half-probabilistic lemma plays a crucial role.\par

Recall that $\sigma_n=\delta_1+\cdots+\delta_n\,$.\par\goodbreak

\begin{lemme}\label{lemme 2.1}
Let $s\geq 2$ and $M$ be integers. Set 
$$\Omega_s(M)=\{\omega\in \Omega\,;\ \Lambda(\omega)\cap [M,+\infty[ 
\ \hbox{contains at least a relation of length s}\,\}\,.$$
Then:
$$\prob\big( \Omega_s(M)\big) 
\leq \frac{s\,2^{s-2}}{ (s-2)!} \sum_{j>M} \delta_j^2\sigma_j^{s-2}
\leq \frac{(4\e)^s}{s^s}\sum_{j>M}\delta_j^2\sigma_j^{s-2}.$$
\end{lemme}

The important fact in this lemma is the presence of the exponent $2$ in the factor $\delta_j^2$ and of the factorial 
in the denominator.\par\medskip

\noindent{\bf Proof.}  We thank the referee for suggesting the following proof.\par
We have $\Omega_s(M)=\bigcup_{l\geq M+s-1} \Delta_l$, where $\Delta_l=\Delta(l,M,s)$ is defined by:
\begin{align*}
\Delta_l 
= \{\omega\,;\  & \Lambda(\omega)\cap [M,+\infty[ \\
& \hbox{ contains at least a relation of length $s$, with largest term $l$}\,\}.
\end{align*}
In other words, $\omega\in\Delta_l$ if and only if $\Lambda(\omega)\cap[M,+\infty[$ has at least a relation of 
length $s$ which contains $l$ and which is contained in $\{M,\ldots,l\}$.\par
We clearly have:
$$\Delta_l\subseteq \bigcup_{(i_1,\ldots,i_{s-1})} \Delta(l,i_1,\ldots,i_{s-1}),$$
where
$$\Delta(l,i_1,\ldots,i_{s-1})
=\{\omega\,;\ \eps_{i_1}(\omega)=\cdots=\eps_{i_{s-1}}(\omega)=\eps_l(\omega)=1\},$$ 
and where $(i_1,\ldots,i_{s-1})$ runs over the $(s-1)$-tuples of integers such that:\par
($\ast$) $M\leq i_1<\cdots<i_{s-1}<l$,\par
($\ast\ast$) $\theta_1i_1+\theta_2i_2+\cdots+\theta_{s-1}i_{s-1} + \theta_s l=0$,
\hskip 3mm $\theta_1,\ldots,\theta_s\in\{-1,+1\}$.\par

Observe that $\delta_{i_{s-1}}\leq (s-1)\delta_l$ for such $(s-1)$-tuples. In fact, it follows from ($\ast\ast$) that
$l\leq i_1+\cdots+i_{s-1}\leq (s-1)i_{s-1}$, so
$$\delta_{i_{s-1}}= \frac{\alpha_{i_{s-1}}}{i_{s-1}}
\leq \frac{\alpha_l}{i_{s-1}}= \frac{\alpha_l}{l} \frac{l}{i_{s-1}}
\leq (s-1)\delta_l.$$\par

Observe also that, when $i_1,\ldots, i_{s-2}$ are fixed, $i_{s-1}=\pm l\pm i_{s-2}\pm\cdots\pm i_1$ can take at 
most $2^{s-1}$ values, so that
\begin{align*}
\prob(\Delta_l)&\leq \sum\prob\big(\Delta(l,i_1,\ldots,i_{s-1})\big)
=\sum\delta_{i_1}\ldots\delta_{i_{s-1}}\delta_l\cr
&\leq (s-1)2^{s-1}\delta_l^2\sum_{M\leq i_1<\cdots<i_{s-2}\leq l-1}
\delta_{i_1}\ldots\delta_{i_{s-2}} \\
& \leq (s-1)2^{s-1}\delta_l^2 \frac{(\delta_M+\cdots+\delta_{l-1})^{s-2}}{(s-2)!}
\end{align*}
by the multinomial formula.\par
Therefore, noting that $s\geq 2$, we have:
\begin{equation}
\prob\big(\Omega_s(M)\big)
\leq\sum_{l\geq M+s-1}(s-1)2^{s-1}\delta_l^2 \frac{\sigma_{l-1}^2}{(s-2)!} 
\leq \sum_{j>M} \frac{(s-1)2^{s-1}}{(s-2)!}\delta_j^2\sigma_j^{s-2}.\tag*{$\square$}
\end{equation}
\medskip

The following theorem is the main result of the paper. It states that subsets of integers can, in several ways, be very 
close to Sidon sets, but in the same time be rather large.\par

\begin{theoreme}\label{theo 2.2}
There exist sets $\Lambda$ of integers which are:\par
{\rm (1)} $p$-Sidon for all $p>1$, $\Lambda(q)$ for all $q<+\infty$, and {\rm CUC}, but which are also\par
{\rm (2)} uniformly distributed; in particular, they are dense in the Bohr group, and ${\cal C}_\Lambda$ 
contains $c_0$, so $\Lambda$ is not a Rosenthal set.
\end{theoreme}

\noindent{\bf Proof.} We use selectors $\eps_k$ of mean 
$$\delta_k=c\, \frac{\log\log k}{k} \qquad (k\geq3),$$ 
where $c$ is a constant to be specified latter. Since this constant plays no role in the beginning of the proof, 
for convenience, we first assume that $c=1$.
\par
The last assertion follows at once from Bourgain's and Lust-Piquard's theorems. The rest of the proof depends on the 
following lemma, where we set $\Lambda_n=\Lambda\cap[1,n]$ and 
$\Lambda'_n=\Lambda\cap[M_n,M_{n+1}[$.

\begin{lemme}\label{lemme 2.3}
If $M_n=n^n$, one has the following properties, where $C_0$ denotes a numerical constant:\par
{\rm (1)} $\sum_{n\geq 3}\prob\big(\Omega_n(M_n)\big)<+\infty.$\par
{\rm (2)} Almost surely $\vert \Lambda_{M_n}\vert\leq C_0\,n(\log n)^2$ for $n$ large enough.\par
{\rm (3)} Almost surely $\vert \Lambda'_n\vert\leq C_0(\log n)^2$ for $n$ large enough.
\end{lemme}

\noindent{\bf Proof of Lemma~\ref{lemme 2.3}.} Note first that
$$\sigma_n=\sum_{3\leq k\leq n} \frac{\log\log k}{k}
\leq (\log\log n)\sum_{3\leq k\leq n} \frac{1}{k}
\leq (\log\log n)\int_1^n \frac{dt}{t}=\log n\,\log\log n.$$\par

Now, take $n\geq 64$, and use Lemma~\ref{lemme 2.1} to obtain (setting $C=4\e$):
\begin{align*}
\prob\big(\Omega_n(M_n)\big)
& \leq  \frac{C^n}{n^n} \sum_{j>M_n} \Big( \frac{\log\log j}{j}\Big)^2 (\log j\,\log\log j)^{n-2} \\
& \leq \frac{C^n}{n^n}\sum_{j>M_n} \frac{(\log j\,\log\log j)^n}{j^2}\,\cdot
\end{align*}

But, for fixed $n$, the function 
$$\frac{u}{v} = \frac{(\log x\,.\log\log x)^n}{x^2}$$ 
decreases on $[M_n,+\infty[$. Indeed, we have to check that $u'(x)v(x)\leq u(x)v'(x)$ on this interval, {\it i.e.} that 
$$nx(1+\log\log x)(\log x\,\log\log x)^{n-1}\leq 2x(\log x\,\log\log x)^n$$
or, equivalently, that
$$n(1+\log\log x)\leq 2\log x\,\log\log x\,.$$
Now, if $x\geq n^n$, we see that 
$$n(1+\log\log x)\leq 2n\log\log x\leq 2n\log n\,\log\log x\leq 2\log x\,\log\log x\,.$$\par

Therefore,
$$\prob\big(\Omega_n(M_n)\big) 
\leq \frac{C^n}{n^n}\int_{M_n}^{+\infty}
\frac{(\log t\,\log\log t)^n}{t^2}\,dt\,.$$\par

Setting 
$$ f_n(t)= \frac{(\log t\,\log\log t)^n}{t^2} \qquad \text{and} \qquad  
I_n=\int_{M_n}^{+\infty} f_n(t)\,dt, $$
we have, by summation by parts:
$$I_n= \frac{(\log M_n\,\log\log M_n)^n}{M_n} 
+\int_{M_n}^{+\infty} nf_n(t)\Big(\frac{1}{\log t}+ \frac{1}{\log t\,\log\log t}\Big)\,dt\,.$$
Since the function in the integrand is less than 
$$\frac{2n}{\log t} f_n(t)\leq \frac{2n}{n\log n} f_n(t)\leq \frac{1}{2}f_n(t)$$
(recall that $n\geq 64$), this gives:
$$I_n\leq \frac{(\log M_n\,\log\log M_n)^n}{M_n}+ \frac{1}{2}I_n\,,$$
so:
\begin{align*}
\prob\big(\Omega_n(M_n)\big)\leq \frac{C^n}{n^n}I_n
&\leq 2 \frac{C^n}{n^n} \frac{(\log M_n\,\log\log M_n)^n}{M_n}\cr
&\leq 2 \frac{C^n}{n^n} \frac{(n\log n.2\log n)^n}{M_n}=
2 \frac{(2C)^n(\log n)^{2n}}{n^n}\,,
\end{align*}
which proves (1).\par

To prove (2), first note that 
$$ \sigma_n\geq \sum_{16\leq k\leq n} \frac{1}{k} \geq \frac{1}{2}\log n \qquad \text{for } n\geq 256.$$ 
Now, using Lemma~\ref{lemme 1.3} with $X_k=\eps_k-\delta_k$, we obtain:
\begin{align*}
\prob\big(\big\vert\,\vert\Lambda_n\vert-\sigma_n\big\vert\geq \frac{\sigma_n}{2}\big) 
& =\prob\Big(\Big\vert \sum_{3\leq k\leq n} X_k\Big\vert \geq \frac{\sigma_n}{2}\Big) \\
& \leq 4\exp\Big(- \frac{\sigma_n}{32}\Big) \leq 4\exp\Big(- \frac{\log n}{64}\Big).
\end{align*}
In particular, since $M_n=n^n$,
$$\prob\Big(\big\vert\,\vert \Lambda_{M_n}\vert - \sigma_{M_n}\big\vert \geq \frac{\sigma_{M_n}}{2}\Big)
\leq 4\exp\Big(- \frac{n\log n}{64}\Big)\,,$$
and the Borel-Cantelli lemma shows that, almost surely,
$$\big\vert\,\vert \Lambda_{M_n}\vert - \sigma_{M_n}\big\vert < \sigma_{M_n}/2$$ 
for $n$ large enough (depending on $\omega$). Thus:
$$\vert \Lambda_{M_n}\vert\leq 2\sigma_{M_n}\leq 2\log M_n\,\log\log M_n \leq C_0\,n(\log n)^2\,,$$
for some numerical constant $C_0$, and this gives (2).\par

The proof of (3) goes the same way. Set:
$$\sigma'_n =\sum_{M_n\leq k<M_{n+1}} \frac{\log\log k}{k}$$
and observe that (here, and in the remainder of the paper, the sign $\sim$ between two functions will mean that 
these two functions are equivalent up to a constant factor):
$$\sigma'_n\sim\log n \sum_{M_n \leq k <M_{n+1}} \frac{1}{k} 
\sim \log n\,\log \frac{M_{n+1}}{M_n} \sim (\log n)^2\,,$$
so that:
$$C_0^{-1}(\log n)^2\leq \sigma'_n\leq C_0(\log n)^2$$
for some numerical constant $C_0$.\par
Then, using again Lemma~\ref{lemme 1.3}, we get:
$$\prob\big(\big\vert\,\vert \Lambda'_n\vert-\sigma'_n\big\vert\geq \sigma'_n/2\big) 
\leq 4\exp\Big(- \frac{\sigma'_n}{32}\Big) 
\leq 4\exp\Big(- \frac{(\log n)^2}{32 C_0}\Big)\, ;$$
so the Borel-Cantelli lemma shows that, almost surely,
$$\vert \Lambda'_n\vert\leq 2\sigma'_n\leq 2C_0(\log n)^2$$
for $n$ large enough, which gives (3), provided we enlarge $C_0$, and completes the proof of 
Lemma~\ref{lemme 2.3}.\hfill\qed
\bigskip

We now conclude the proof of Theorem~\ref{theo 2.2}.\par
We first choose the constant $c$ in order that, not only 
$\sum_{n\geq 3}\prob\big(\Omega_n(M_n)\big)<+\infty$, but 
$\sum_{n\geq 3}\prob\big(\Omega_n(M_n)\big) < 1$. So, using Lemma~\ref{lemme 2.3}, we can find 
$\Omega_0\subseteq \Omega$ such that $\prob(\Omega_0)>0$ and with the property:
\begin{equation}\label{eq 4}
\begin{split}
&\text{If $\omega\in\Omega_0$, then $\omega\not\in \bigcup\limits_{n\geq 3}\Omega_n(M_n)$}.\quad 
\text{There exists $n_0=n_0(\omega)$ such that}\\ 
& \text{$\vert \Lambda_{2M_n}\vert\leq C_0\,n(\log n)^2$ 
and $\vert \Lambda'_n\vert\leq C_0(\log n)^2\leq n$ for $n>n_0$.}
\end{split}
\end{equation}
Indeed, an inspection of the proof of Lemma~\ref{lemme 2.3} shows that we also have, almost surely, 
$\vert \Lambda_{2M_n}\vert \leq 2\sigma_{2M_n}\leq C_0\,n(\log n)^2$ for $n$ large enough, and this gives 
\eqref{eq 4}. We have the following consequences, where $\omega\in \Omega_0$, and 
$\Lambda=\Lambda(\omega)$:
\begin{equation}\label{eq 5}
\Lambda\cap[M_n,+\infty[ \text{ contains no relation of length $\leq n$.}
\end{equation}
For, if $\Lambda\cap[M_n,+\infty[$ were to contain a relation $R$ of support $S$ with $\vert S\vert=s\leq n$, then 
necessarily $s\geq 3$, $S\subseteq \Lambda\cap[M_s,+\infty[$ and $\omega\in\Omega_s(M_s)$, which is not 
the case by \eqref{eq 4}. Now, \eqref{eq 4} and \eqref{eq 5} imply that, for $n>n_0$, $\Lambda'_n$ is 
quasi-independent, and so is a Sidon set with bounded constant. So, we get that $\Lambda$ is CUC and 
$\Lambda(q)$ for all $q<+\infty$ using Propositions~\ref{prop 1.1} and \ref{prop  1.2}, provided we notice that:
$$ \frac{M_{n+1}}{M_n}= \frac{(n+1)^{n+1}}{n^n}\geq n+1\geq 2.$$\par

To end the proof, we first show that $\Lambda$ is $p$-Rider, for every $p>1$, and then, using 
Proposition~\ref{prop 1.6} (c), prove that it is $q$-Sidon for every $q>1$.\par\smallskip

So, fix $p\in]\,1,2\,[$, set $\eps=2/p-1\in]\,0,1\,[$, and take $\omega\in \Omega_0$ and 
$n_1=n_1(\eps,\omega)\geq 2n_0(\omega)$ such that $C_0\,n(\log n)^2\leq n^{1/\eps}/2$ and 
$n^{1/\eps}/2\geq n$ for $n\geq n_1$.\par

Let $A\subseteq\Lambda$ be a finite subset, with $\vert A\vert^\eps>n_1$. Set $n=[\vert A\vert^\eps]$, where 
$[\ ]$ stands for integer part, so that $n\geq n_1$ and $\vert A\vert\geq n^{1/\eps}$. Observe that:
\begin{align*}
\big\vert A\cap[M_n,+\infty[\big\vert 
& \geq \vert A\vert - \vert A\cap[1,M_n]\vert \geq \vert A\vert-\vert \Lambda_{M_n}\vert \\
& \geq n^{1/\eps}-C_0\,n(\log n)^2\geq n^{1/\eps}/2\geq n\,,
\end{align*}
and select $B\subseteq A\cap[M_n,+\infty[$ with $\vert B\vert=n$. It follows from \eqref{eq 5} that $B$ is 
quasi-independent, and  $\vert B\vert=n \geq\frac{1}{2}\vert A\vert ^\eps$.\par 
If now $A$ is a subset of $\Lambda$ with $1\leq \vert A\vert\leq n_1$, simply take for $B$ a singleton from 
$A$. Then $B$ is quasi-independent, and $\vert B\vert=1\geq n_1^{-1}\vert A\vert^\eps$.\par
The criterion of Proposition~\ref{prop 1.6} (a) is verified with $\delta=n_1^{-1}$. Therefore $\Lambda$ is 
$p$-Rider.\par

We shall verify that we are in position to apply part (b) of Proposition~\ref{prop 1.6}.\par
Take $p\in]\,1,2\,[$ and $1/p<\eps<1$. Take $\omega\in\Omega_0$ and 
$n_1=n_1(\eps,\omega)\geq 2n_0(\omega)$ such that $C_0\,n(\log n)^2\leq n^{1/\eps}/2$ and 
$n^{1/\eps}/2\geq n$ for $n\geq n_1$.\par
Let $A\subseteq \Lambda$ be a finite subset with $\vert A\vert^\eps>n_1$. Set $n=[\vert A\vert^\eps]$, where 
$[\ ]$ stands for integer part, so that $n\geq n_1$ and $\vert A\vert\geq n^{1/\eps}.$ Observe that:
\begin{align*}
\vert A\cap[2M_n,+\infty[\vert 
&\geq\vert A\vert-\vert A\cap[1, 2M_n]\vert
\geq\vert A\vert-\vert \Lambda_{2M_n}\vert \\
& \geq n^{1/\eps}-C_0\,n(\log n)^2 \geq n^{1/\eps}/2\geq n
\end{align*}
in view of \eqref{eq 4}. We can thus select $B\subseteq A\cap[2M_n,+\infty[$ with $\vert B\vert=n-1$ and have:
\begin{displaymath}
\text{If $k\in \Lambda\cap [M_n,+\infty[$, then $B\cup\{k\}$ is quasi-independent.}
\end{displaymath}
Indeed, $B\cup\{k\}$ is a set of cardinality less than $n$ contained in $\Lambda\cap[M_n,+\infty[$, and is 
automatically quasi-independent, from \eqref{eq 5}.\par
\medskip

We show now that $B$ is 8-pseudo-complemented in $\Lambda$.\par
Put $\nu=\delta_0-V_{M_n}$, where $\delta_0$ is the Dirac point mass at $0$, and $V_{M_n}$ the de la 
Vall\'ee-Poussin kernel of order $M_n$. Consider the Riesz product  $R=\prod\limits_{k\in B}(1+\Re e_k)$, and 
set $\mu=2\nu\ast R$. We claim that:
$$\Vert \mu\Vert\leq 8\,; \qquad \hat\mu\geq 1 \quad \text{on } B\,; \qquad \hat\mu=0 \quad 
\text{on } \Lambda\setminus B.$$
Indeed, $\Vert\nu\Vert\leq 4$ and $B$ is quasi-independent, so the Riesz product $R$ verifies 
$\Vert R\Vert=\hat R(0)=1$. Therefore $\Vert\mu\Vert \leq 8$.\par
Take $l\in B$. Then $l> 2M_n$ and $\hat\nu(l)=1$. As $\hat R(l)\geq 1/2$, we have $\hat\mu(l)\geq1$.\par
If $A\subseteq\Lambda$ and $\vert A\vert^\eps\leq n_1$, any singleton $B$ of $A$ is quasi-independent, 
$1$-complemented in $\Lambda$, and $\vert B\vert \geq n_1^{-1}\vert A\vert^\eps$.\par
\medskip

We have thus verified the hypothesis of part (b) of Proposition~\ref{prop 1.6}, and so $\Lambda$ is $q$-Sidon 
for any $q>1/\eps$. In particular, it is $p$-Sidon, and this ends the proof of Theorem~\ref{theo 2.2}.\hfill\qed\par
\bigskip

\noindent{\bf Remark 1.} The proof shows that we can actually extract from $A$, for every $\alpha>0$, a 
quasi-independent set $B$ such that $\vert B\vert\geq \delta \vert A\vert/(\log\vert A\vert)^{2+\alpha}$. 
Moreover, a slight modification leads to sets even closer to Sidon sets. 
\par

\begin{proposition}\label{prop 2.4}
Let $\alpha>1$ and $\varphi_\alpha$ be the Orlicz function $x\mapsto x\big(\log(1+x)\big)^\alpha$. Then, there 
exists a set $\Lambda$ as in Theorem~\ref{theo 2.2}, and moreover such that 
$\hat f\in \ell_{\varphi_\alpha,\infty}$ for every $f\in{\cal C}_\Lambda$.
\end{proposition}

Recall that $\ell_{\varphi_\alpha,\infty}$ is the weak Orlicz-Lorentz space of sequences $(a_n)_n$ such that 
$\sup\limits_n \varphi_\alpha^{-1}(n)a_n^\ast<+\infty$, where $(a_n^\ast)_n$ is the non-increasing 
rearrangement of $(\vert a_n\vert)_n$. Therefore, another way to phrase the proposition is, setting 
$a_n=\hat f(n)$:
$$a_n^\ast\leq C_\alpha\Vert f\Vert_\infty(\log n)^\alpha/n \quad \text{for every } f\in{\cal C}_\Lambda.$$
The proof just consists in changing $M_n$. We take $M_n=\big[e^{n(\log\log n)^2}\big]$, where $[\;]$ stands 
for the integer part. We still have $\sum\limits_n\prob\big(\Omega_n(M_n)\big)<+\infty$, since 
$$\prob\big(\Omega_n(M_n)\big) 
\leq \frac{2\,C^n}{n^n} \frac{(\log M_n\,\log\log M_n)^n}{M_n}
\leq  \exp\big(-n\,(\log\log n)^2/2\big)$$
for $n$ large enough. Arguing as previously, we get for every finite subset $A$ of $\Lambda$, a quasi-independent 
subset $B$ of $A$ such that $\vert B\vert\geq\delta\vert A\vert/(\log\vert A\vert)^\alpha$, and such that 
${\cal C}_B$ is uniformly pseudo-complemented in ${\cal C}_A$. As in the proof of Proposition~\ref{prop 1.6}, we 
obtain
$$\vert\,\{\vert\hat f\vert>t\}\,\vert \leq C\,\varphi_\alpha\Big(\frac{\Vert f\Vert_\infty}{t}\Big),$$
which gives the result (arguing as in \cite{Lefevre2} for instance).\par
\medskip

We cannot eliminate a logarithmic factor, and replace $\alpha>1$ by $\alpha>0$ because, due to Bourgain's 
criterion, we have to assume that $\sigma_n/\log n$ goes to infinity in order that ${\cal C}_\Lambda$ contains 
$c_0$. However, for each $\alpha>0$, there do exist non-Sidon sets 
$\Lambda$ for which $\hat f\in \ell_{\varphi_\alpha}$ when $f\in {\cal C}_\Lambda$ (as can be seen from 
\cite{Blei-K}, p. 69).\par
The set $\Lambda$ is, in some sense, very close to be Sidon, whereas ${\cal C}_\Lambda$ contains $c_0$. 
However, it cannot be too close without being Sidon because if $\hat f\in \ell_{1,\infty}$, the Lorentz space 
weak-$\ell_1$, for every $f\in {\cal C}_\Lambda$, then $\Lambda$ is Sidon. In fact, this condition implies an 
inequality of the type:
\begin{equation}\label{etoile}
\vert\{\vert\hat f\vert\geq t\}\vert \leq \frac{C}{t}\Vert f\Vert_\infty
\tag*{$(\ast)$}
\end{equation}
for every $f\in {\cal C}_\Lambda$. Let now $A$ be a finite subset of $\Lambda$, and 
$f=\sum\limits_{n\in A}e_n$ and $f_\omega=\sum\limits_{n\in A}r_n(\omega)e_n$, where $r_n$, $n\geq 1$ 
are the Rademacher functions. Then, inequality \eqref{etoile} applied with $t=1$ gives 
$\Vert f_\omega\Vert_\infty\geq (1/C)\vert A\vert$. Integrating in $\omega$ gives 
$\lpnorm f\rpnorm\geq (1/C)\vert A\vert$, from which follows, by a result of G. Pisier 
(\cite{Pisier48}, Th\'eor\`eme 2.3 (vi)), that $\Lambda$ is a Sidon set.\par
\bigskip

\noindent{\bf Remark 2.} If one takes selectors of mean $\delta_n$ such that $n\delta_n$ is bounded, the 
corresponding random set $\Lambda(\omega)$ is almost surely a Sidon set. This is a well-known result of 
Y. Katznelson and P. Malliavin (\cite{Katz-Mall}, or \cite{Katz}), and Lemma~\ref{lemme 2.1} gives another proof 
of this fact. It suffices to take $M_n=A^n$, where $A$ is a given integer, large enough to have 
$\sum_{n=1}^{+\infty}\prob\big(\Omega_n(M_n)\big) <1$. Then, with positive probability 
$\Lambda(\omega)\cap[M_n,+\infty[$ contains no relation of length $\leq n$, whereas 
$\vert \Lambda(\omega)\cap[1,M_n]\vert\leq C n$. Hence, for every finite subset $A$ of $\Lambda(\omega)$, we 
can find a quasi-independent subset $B\subseteq A$ such that $\vert B\vert\geq\delta\vert A\vert$, for some 
fixed $\delta=\delta(\omega)$. It follows from Pisier's characterization (\cite{Pisier48}, Th. 2.3 (iv)) that, with 
positive probability, and hence almost surely by Kolmogorov's $0-1$ law, $\Lambda(\omega)$ is a Sidon set.\par
\smallskip

As is now well-known, Sidon sets are characterized by various properties (successively weaker) of the Banach space 
${\cal C}_\Lambda$: 
$\Lambda$ is a Sidon set {\it iff} ${\cal C}_\Lambda$ is isomorphic to $\ell_1$ (\cite{Varo}), {\it iff} 
${\cal C}_\Lambda$ has cotype 2 (\cite{Kwap-P}, Th. 3.1, \cite{Pisier1}), and {\it iff} ${\cal C}_\Lambda$ has a 
finite cotype (\cite{Bourgain-Milman}). This later property can be expressed by saying that ${\cal C}_\Lambda$ 
does not contain $\ell_\infty^n$ uniformly. So, deterministically, one has the 
dichotomy:\par
(a) either $\Lambda$ is a Sidon set, and so ${\cal C}_\Lambda$ is isomorphic to $\ell_1$;\par
(b) or ${\cal C}_\Lambda$ contains $\ell_\infty^n$ uniformly.\par

The probabilistic dichotomy is stronger: taking selectors of mean $\delta_1,\delta_2,\ldots\,$, with $(\delta_n)_n$ 
decreasing, one has:\par
(a) either almost surely $\Lambda$ is a Sidon set (if $n\delta_n$ is bounded);\par
(b) or almost surely ${\cal C}_\Lambda$ contains $c_0$ (if $n\delta_n$ is not bounded), and $\Lambda$ is even 
uniformly distributed.\par
\medskip

Y. Katznelson (\cite{Katz}) already noticed such a ``dichotomy'': he showed that (under a different choice of 
selectors from ours) either almost surely $\Lambda$ is a Sidon set, or almost surely $\Lambda$ is dense in the 
Bohr group. However, this is perhaps not a true dichotomy since it is a well-known open problem whether 
there can exist Sidon sets dense in the Bohr group (see \cite{Myriam}, question 2, p. 14; it is stated for the Bohr 
group of $\R$, but also makes sense for the Bohr group of $\Z$).\par

The dichotomy stated here strengthens Katznelson's result since every uniformly distributed set is dense in the 
Bohr group (see \cite{Blum}, Theorem 1); indeed, saying that $\Lambda=\{\lambda_1,\lambda_2,\ldots\}$ is 
uniformly distributed means that the measures $\mu_N=1/N\sum\limits_{n=1}^{+\infty}\delta_{\lambda_n}$ 
($\delta_{\lambda_n}$ is there the Dirac measure at the point $\lambda_n$) converge weak-star to the Haar 
measure $\mu$ of the Bohr group $b\Z$; but these measures are carried by $\Lambda$, so the closed support of 
$\mu$ is contained in the Bohr closure of $\Lambda\,;$ since the Haar measure is continuous, we get that this 
closure is the whole Bohr group.\par
\bigskip

\noindent{\bf Remark 3.} The random sets $\Lambda$ that we construct have an asymptotical 
quasi-independence: $\Lambda\cap[M_n,+\infty[$ contains no relation of length $\leq n$. This is reminiscent of 
the following result of J. Bourgain (\cite{Bourgain7}): if $\Lambda$ is a Sidon set and $n\in\N^\ast$, there exists 
$l_n=l(\Lambda,n)$ such that $\Lambda$ can be decomposed in $l_n$ sets $\Lambda_1,\ldots,\Lambda_{l_n}$, 
each of which contains no relation of length $\leq n$.\par
\bigskip

We now investigate what happens when we let $p$ increase away from $1$.We get several different results, 
and $p=4/3$ seems to play a special role.\par
\bigskip

We first state two very similar results.\par

\begin{theoreme}\label{theo 2.5}
For every $1<p < 4/3$, there exists a set $\Lambda$ of integers which is:\par
{\rm (1)} uniformly distributed (so $\Lambda$ is dense in the Bohr group, ${\cal C}_\Lambda$ contains $c_0$, and 
$\Lambda$ is not a Rosenthal set), and which is:\par
{\rm (2)} $\Lambda(q)$ for all $q<+\infty$, a {\rm CUC}-set, and moreover is:\par
\qquad {\rm (a)} $p$-Rider, but not $q$-Rider for $q<p$\par
\qquad {\rm (b)} $q$-Sidon for all $q>p/(2-p)$.
\end{theoreme}

\begin{theoreme}\label{theo 2.6}
Same as Theorem~\ref{theo 2.5}, except that, instead of {\rm (a)}, $\Lambda$ is:\par
{\rm (a')} $q$-Rider for every $q>p$, but is not  $p$-Rider.
\end{theoreme}

\noindent{\bf Remark.} After this paper was completed, P. Lef\`evre and the third-named author proved  
(\cite{Lefevre-Luis}) that every $p$-Rider set with $p<4/3$ is a $q$-Sidon set, for all $q>p/(2-p)$. A weaker, 
unpublished, result, due to J. Bourgain, is quoted in \cite{Myriam}, p. 41. Hence condition {\rm (b)} always 
follows  from condition {\rm (a)}, and is not specific to the construction. We do not know whether this gap 
between $p$ and $p/(2-p)$ follows only from technical reasons. For $p>1$, whether every $p$-Rider set is 
actually $p$-Sidon is an open question.\par
In Theorem~\ref{theo 2.5}, we obtain sets which are $p$-Rider but not $q$-Rider for $q<p$. We do not know if 
these sets are $p$-Sidon, so exactly $p$-Sidon, in the terminology of R. Blei. He constructed such sets using 
fractional products (\cite{Blei2}, \cite{Blei3}). We may call the sets in Theorem~\ref{theo 2.5} ``{\it exactly 
$p$-Rider sets}''. The sets appearing in Theorem~\ref{theo 2.6} are of a different kind. We may call them 
``{\it exactly $p^+$-Rider sets}''. Such sets were also obtained in \cite{Blei2}, Corol. 1.7 d), where they were 
called ``exactly non-$p$-Sidon'', and were called ``asymptotic $p$-Sidon'' in \cite{Blei-K}.\par
\bigskip

\noindent{\bf Proof.} It is similar to that of Theorem~\ref{theo 2.2}, so we shall be 
very sketchy.\par
Let $\alpha=2(p-1)/(2-p)\in]\,0,1\,[$.\par
For Theorem~\ref{theo 2.5}, we use selectors $\eps_k$ of mean 
$$ \delta_k=c\, \frac{(\log k)^\alpha}{k(\log\log k)^{\alpha+1}} \qquad \text{for } k\geq 4.$$ 
As in Lemma~\ref{lemme 2.3}, we have, with $M_n=n^n$, 
$\sum_{n\geq 1}\prob\big(\Omega_n(M_n)\big)<+\infty$, and almost surely
$C_0 n^{\alpha+1}\leq \vert\Lambda_{M_n}\vert \leq C_1\,n^{\alpha+1}$ and 
$\vert\Lambda'_n\vert\leq C_1\,n^\alpha$ for $n$ large enough.\par

For Theorem~\ref{theo 2.6}, we increase the means $\delta_k$ slightly, replacing them by 
\begin{equation}
\delta_k=c\, \frac{(\log k)^\alpha\log\log k}{k}\,\cdot \tag*{$\square$}
\end{equation}

\noindent{\bf Remark.} In order to prove our theorems, we used selectors with various means. They are smaller in 
Theorem~\ref{theo 2.2} than in Theorem~\ref{theo 2.5}, for instance. We remark that selectors $(\eps_k)_k$ 
of mean $\delta_k$ with $\delta_k\leq \delta'_k$ may be achieved as the product of two independent sequences of 
selectors $(\eps'_k)_{k}$ and $(\eps''_k)_{k}$ of mean $\delta'_k$ and $\delta''_k=\delta_k/\delta'_k$. 
It follows that, for example, the sets in Theorem~\ref{theo 2.2} may be constructed inside the respective sets of 
Theorem~\ref{theo 2.5}.\par
\bigskip

In Theorem~\ref{theo 2.5}, the proof that $\Lambda$ was {\it CUC} or $\Lambda(q)$ was based on the fact that 
$\vert\Lambda_n'\vert\subseteq\Lambda\cap[M_n,+\infty[$ is quasi-independent. For $\alpha\geq 1$ ({\it i.e.} 
$p\geq 4/3$), we no longer have $\vert\Lambda_n'\vert\leq n$, and therefore, {\it a priori}, must give up these 
properties. However, we can use another extraction procedure. This procedure was first introduced by J. Bourgain 
(\cite{Bourgain8}); later, a clear statement was given in \cite{Luis2}, \S~III.2. Since this last reference is hardly 
available, we prefer to give a self-contained proof.\par 

The corresponding set $\Lambda(\omega)$ of integers that we shall obtain in this manner satisfies  
$\vert\Lambda(\omega)\cap[2^n,2^{n+1}[\vert\sim n\sim \log 2^n$, which is the limiting  condition of mesh 
(on arithmetic progressions) for Sidon sets. This size is in some sense the largest possible if we want 
to obtain a set $\Lambda$ with blocks having a uniformly bounded Sidon constant.\par\goodbreak

\begin{theoreme}\label{theo 2.7}
There exists a set $\Lambda$ of integers which is uniformly distributed and contains a subset 
$E\subseteq\N^\ast$ which is:
\begin{itemize}
\item[{\rm (1)}] $4/3$-Rider, and not $q$-Rider for $q<4/3$; a {\rm CUC}-set; a $\Lambda(q)$-set for all 
$q<+\infty$ (more precisely, for all $q>2$, we have: $\Vert f\Vert_q\leq Cq^2\Vert f\Vert_2$ for all 
$f\in{\cal P}_E$, where $C>0$ is a numerical constant), and nevertheless, 
\item[{\rm (2)}] has positive upper density in $\Lambda$, so, in particular, ${\cal C}_E$ contains $c_0$, 
and $E$ is not a Rosenthal set.
\end{itemize}
\end{theoreme}

Let $A$ be a finite subset of integers. For the proof, it will be convenient to define:
$$\psi_A=\sup_{p\geq2} \frac{\Vert e_A\Vert_p}{\sqrt p}, 
\qquad \text{where } e_A=\sum\limits_{k\in A}e_k\,.$$
\par

We need the following simple estimate of $\psi_A$.

\begin{lemme}\label{lemme 2.8}
Let $I=[a+1,a+N]$ be an interval of integers of length $N$, $N\geq 3$. Then:
$$\psi_I\leq \frac{N}{\sqrt{2\log N}} \, \cdot$$
\end{lemme}

\noindent{\bf Proof.} For $p\geq 2$,  $\vert e_I\vert^p\leq N^{p-2}\vert e_I\vert^2$, so  
$\int\vert e_I\vert^p\,dm \leq N^{p-2}\int\vert e_I\vert^2\,dm=N^{p-1}$ and 
$\Vert e_I\Vert_p / \sqrt p \leq N^{1-1/p} / \sqrt p$. Optimizing gives $p=2\log N$ ($\geq 2$), and 
the lemma.\hfill\qed
\bigskip

This estimate is essentially optimal. Indeed, it is well-known that $\psi_I$ is uniformly equivalent to 
$\theta=\Vert e_I\Vert_\Psi$ ($\Vert\;\Vert_\Psi$ being the norm associated to the Orlicz function 
$\Psi(x)=\e^{x^2}-1$). But, for some constant $\gamma$, $\vert e_I(t)\vert\geq \gamma N$ for $t$ in 
an interval $J$ of length $\geq \gamma N^{-1}$ around $0$, so one has:
$$2\geq\int_J\exp\Big( \frac{\vert e_I\vert^2}{\theta^2}\Big)\,dm 
\geq \gamma N^{-1}\exp\Big(\frac{\gamma^2 N^2}{\theta^2}\Big)\,,$$
whence $\theta\geq \gamma^{-1}N / \sqrt{\log 2\gamma^{-1}N}$.
\par\smallskip

We now use selectors $\eps_k$ of mean $\delta_k=c\,n/2^n$ for $2^n\leq k<2^{n+1}$, where $c>0$ is a 
given constant.\par
Set
$$I_n=[2^n,2^{n+1}[,\ n\geq2\,;\hskip 5mm\delta_k=c\, \frac{n}{2^n}\hbox{ if }k\in I_n\,.$$
\par

Note that $(\delta_k)_k$ decreases, and $\delta_k$ is of the form $\alpha_k / k$, where $(\alpha_k)_k$ goes to 
$+\infty$.\par

If $\Lambda=\Lambda(\omega)$ is the corresponding set of integers, it will be convenient to set:
$$\Lambda_n=\Lambda\cap I_n\,;\hskip 5mm 
\sigma_n=\esp\vert \Lambda_n\vert=\sum_{k\in I_n}\delta_k=cn\,.$$\par

For this proof, the value of $\psi_{I_n}$ is somewhat large, and requires $c$ be sufficiently small, say $c\leq 1/576$. 
We prefer to follow another route, which could be useful in other contexts, by choosing also a random set in $I_n$ 
for which the $\psi$ constant is small enough. We make the two random choices at the same time. Namely, 
we consider $(\eps'_n)_{n\geq1}$, a second sequence of selectors, independent of $(\eps_n)_{n\geq1}$, with fixed 
mean $\tau$, and set $\Lambda'_n (\omega) = \{k \in \Lambda_n (\omega)\,; \eps_k' (\omega) = 1\}$. In 
short:
$$\Lambda'_n=\{k\in \Lambda_n\,;\ \eps'_k=1\}\,;\hskip 5mm
\Lambda'=\bigcup_{n=1}^{+\infty} \Lambda'_n\,.$$

The following lemma, which is a slight modification of Bourgain's construction in \cite{Bourgain8}, is really the 
heart of the proof.\par

\begin{lemme}\label{lemme 2.9}
Almost surely, for $n$ large enough, one has:\par
{\rm (1)} $(c/2)\,n\leq \vert\Lambda_n\vert\leq (2c)\,n$ and 
$(c\tau/2)\,n\leq \vert \Lambda'_n\vert\leq (2c\tau)\,n$\par
{\rm (2)} $\Lambda'_n$ contains at most relations of length $\leq l_n$, where $l_n=[144\,c^2\tau^2n]$.
\end{lemme}

\noindent{\bf Proof of Lemma~\ref{lemme 2.9}.} We have already seen that: 
$$\prob\Big(\big\vert\, \vert\Lambda_n\vert-\sigma_n\big\vert
\geq \frac{\sigma_n}{2}\Big)\leq\exp\Big(- \frac{\sigma_n}{32}\Big) =
\exp\Big(- \frac{c\,n}{32}\Big)\,,$$ 
so, by the Borel-Cantelli lemma, $\vert \Lambda_n\vert$ is almost surely between $(c/2)\,n$ and $(2c)\,n$ for 
$n$ large enough; and this proves the first half of (1). The second half holds for the same reason, since 
$\Lambda'_n$ corresponds to selectors $\eps_k\eps'_k$ with mean $(c\tau)\,n/2^n$ for $k\in I_n$.\par
The proof of (2) is more elaborate.\par
Fix $n$, and consider the random trigonometric polynomial:
$$F_\omega=\sum_{j=l_n+1}^{\vert I_n\vert}
\sum_{{\mathop{\scriptscriptstyle R\subseteq I_n}\limits_{\scriptscriptstyle \vert R\vert = j}}} 
\prod_{k\in R} \eps_k(\omega)\eps'_k(\omega)\big(e_k+e_{-k}\big)\,.$$
\par

Set:
$$N_n(\omega) = \int_\T F_\omega(t)\,dm(t)\,.$$
Expanding $F_\omega$, we see that:
\begin{align*}
F_\omega(t) 
& =\hskip -1mm\sum_{j=l_n+1}^{\vert I_n\vert} 
\sum_{{\mathop{\scriptscriptstyle R\subseteq I_n}\limits_{\scriptscriptstyle \vert R\vert = j}}}
\ \sum_{\theta_k\in\{-1,+1\}^R}
\prod_{k\in R}\eps_k(\omega)\eps'_k(\omega)e_k^{\theta_k}(t) \cr
&=\hskip -1mm\sum_{j=l_n+1}^{\vert I_n\vert} 
\sum_{{\mathop{\scriptscriptstyle R\subseteq I_n}\limits_{\scriptscriptstyle \vert R\vert = j}}}
\ \sum_{\theta_k\in\{-1,+1\}^R}\e^{it\big(\sum_{k\in R}\theta_k k\big)}.
\end{align*}
\par

The contribution to $N_n(\omega)$ of an exponential of this sum is $0$ if 
$\sum\limits_{k\in R}\theta_k k\not=0$, and is $1$ if $\sum\limits_{k\in R}\theta_k k=0$. Therefore, 
$N_n(\omega)$ is exactly the number of relations of length $>l_n$ in $\Lambda'_n$.\par

We claim that $N_n(\omega)$ is almost surely zero for $n$ large enough. To that effect, we majorize the 
expectation $J$ of $N_n(\omega)$, using Fubini's theorem. Indeed, $J=\int_\T H(t)\,dm(t)$, where:  
$$H(t)=\int_\Omega F_\omega(t)\,d\prob(\omega)=
\sum_{j=l_n+1}^{\vert I_n\vert} 
\sum_{{\mathop{\scriptscriptstyle R\subseteq I_n}\limits_{\scriptscriptstyle \vert R\vert = j}}}
\delta^j\prod_{k\in R}(e_k+e_{-k})\ .$$
and $\delta=c\tau\,n/2^n$. Hence:
$$J=\sum_{j=l_n+1}^{\vert I_n\vert} 
\sum_{{\mathop{\scriptscriptstyle R\subseteq I_n}\limits_{\scriptscriptstyle \vert R\vert = j}}} \delta^j
\int_\T \prod_{k\in R}\big((e_k(t)+e_{-k}(t)\big)\,dm(t)\,.$$
At this stage, it is useful to observe that:
\begin{equation}\label{eq 3}
\begin{split}
\sum_{{\mathop{\scriptscriptstyle R\subseteq I_n}\limits_{\scriptscriptstyle \vert R\vert = j}}}
\int_\T \prod_{k\in R}\big( e_k (t) + e_{-k} (t) \big) 
& \,dm(t)  \cr
& \leq \frac{1}{ j!} \int_\T \Big(\sum_{k\in I_n} \big((e_k(t) + e_{-k}(t)\big)\Big)^j\,dm(t)\,.
\end{split}
\tag*{{\rm (3)}}
\end{equation}
Indeed, when we expand 
$$\Big(\sum_{k\in I_n}\big((e_k(t)+e_{-k}(t)\big)\Big)^j \,,$$ 
each term $\prod_{k\in R}\big(e_k(t)+e_{-k}(t)\big)$ appears $j!$ times, whereas the other terms on the right 
hand side of \ref{eq 3} are positive. It now follows from \ref{eq 3} that:
\begin{align*}
J\leq \sum_{j=l_n+1}^{\vert I_n\vert} \frac{\delta^j}{j!} \int_\T \Big(\sum_{k\in I_n}
&\big((e_k (t) + e_{-k}(t)\big)\Big)^j\,dm(t) \cr
&\leq \sum_{j=l_n+1}^{\vert I_n\vert} \frac{\delta^j}{j!} 2^j \Big\Vert \sum_{k\in I_n} e_k\Big\Vert_j^j
\leq \sum_{j=l_n+1}^{\vert I_n\vert} \frac{2^j\delta^j}{ j!} (\psi_{I_n}\sqrt j)^j.
\end{align*}
Since $j!\geq (j / \e)^j \geq (j / 3)^j$, this gives
$$J\leq \sum_{j=l_n+1}^{+\infty} \Big( \frac{6\delta\psi_{I_n}}{\sqrt j}\Big)^j
\leq \sum_{j=l_n+1}^{+\infty} \Big(\frac{6\delta\psi_{I_n}}{\sqrt {l_n+1}}\Big)^j.$$
Therefore, 
$$J\leq 2^{-l_n} \qquad  \text{if} \quad  \frac{6\delta\psi_{I_n}}{\sqrt{l_n+1}}\leq \frac{1}{2}\, 
\raise 1pt \hbox{,}$$ 
{\it i.e.} if $l_n+1\geq 144 (\delta\psi_{I_n})^2$. But, it follows from Lemma~\ref{lemme 2.8} that:
$$\psi_{I_n}\leq \frac{2^n}{\sqrt{(2\log2)\,n}}\leq \frac{2^n}{\sqrt n}\,\cdot$$
Therefore 
$$ 144(\delta\psi_{I_n})^2\leq 144\Big(c\tau \frac{n}{2^n}\cdot \frac{2^n}{\sqrt n}\Big)^2 
= 144\,c^2\tau^2 n,$$ 
and the choice of $l_n$ just fits to obtain $J \leq 2^{-l_n}$. Of course, we have assumed $n$ large enough to have 
$l_n\geq 1$ in that proof.\par\smallskip

Finally, Markov's inequality implies:
$$\sum_{n\geq2} \prob(N_n\geq1) \leq\sum_{n\geq 2}\esp N_n \leq \sum_{n\geq2} 2^{-l_n} <+\infty\,,$$
and by the Borel-Cantelli lemma, the integer $N_n$ is almost surely zero for $n$ large enough, and that ends the 
proof of Lemma~\ref{lemme 2.9}.\hfill\qed
\bigskip

Now, using Bourgain's Theorem~\ref{theo 1.10} and Lemma~\ref{lemme 2.9}, one can find 
$\Omega_0\subseteq\Omega$ with $\prob(\Omega_0)=1$ such that for $\omega\in\Omega_0$, there exists 
$n_0=n_0(\omega)$ such that $\Lambda=\Lambda(\omega)$ and $\Lambda'=\Lambda'(\omega)$ satisfy:\par
\begin{itemize}
\item[{\rm (4)}] $\Lambda$ and $\Lambda'$ are uniformly distributed\par
\item[{\rm (5)}] $(c/2)\,n\leq\vert\Lambda_n\vert\leq (2c)\,n$ and 
$(c\tau/2)\,n\leq \vert\Lambda'_n\vert\leq (2c\tau)\,n$ for $n>n_0$\par
\item[{\rm (6)}] $\Lambda'_n$ contains at most relations of length less than $\leq l_n=[144c^2\tau^2n]$ for 
$n>n_0$.\par
\end{itemize}

$\Lambda'_n$ is not quite quasi-independent, so we shall modify it slightly. We adjust once and for all $\tau$, 
depending on $c$, such that $144c^2\tau^2\leq c\tau/4$ ({\it e.g.} taking $c\tau=1/576$), so that 
$l_n\leq c\tau\,n/4\leq \vert \Lambda'_n\vert/2$ for $n>n_0$, in view of (5). Select then in $\Lambda'_n$ a 
relation $R$ with support $S_n$ of maximal cardinality. Then $\vert S_n\vert\leq l_n$ from (6), and 
$E_n=\Lambda'_n\setminus S_n$ is quasi-independent. Moreover:
$$\vert E_n\vert=\vert\Lambda'_n\vert-\vert S_n\vert 
\geq\vert\Lambda'_n\vert-l_n\geq\vert\Lambda'_n\vert/2$$
for $n>n_0$. Hence, if we set $E=\bigcup\limits_{n>n_0}E_n$, we have $E_n=E\cap I_n$, and, moreover:\par
\begin{itemize}
\item[{\rm (7)}] $E$ has positive upper density in $\Lambda$\par
\noindent (note that $\Lambda'$ has upper density $\geq\tau/4$ in $\Lambda$ by (5)),\par
\item[{\rm (8)}] $E_n$ is quasi-independent,\par
\item[{\rm (9)}] $\vert E_n\vert\geq (c\tau/4)\,n$,\par
\item[{\rm (10)}] If $A\subseteq E$ is a finite subset, then $A$ contains a quasi-independent subset $B$ with 
$\vert B\vert\geq(1/2)\,\vert A\vert^{1/2}$.\par
\end{itemize}

The last property is proved in the following way. Set $Z=\{n\,;\ A\cap E_n\not=\emptyset\}$ and 
$h=\vert Z\vert$. We distinguish two cases.\par
\smallskip

{\it Case 1: there exists $n\in Z$ such that $\vert A\cap E_n\vert\geq \vert A\vert^{1/2}$.}\par
Then, just take $B=A\cap E_n$ to have a quasi-independent set $B$ such that 
$\vert B\vert\geq \vert A\vert^{1/2}$.\par\smallskip

{\it Case 2: $\vert A\cap E_n\vert<\vert A\vert^{1/2}$ for any $n\in Z$.}\par
Then $h\geq \vert A\vert^{1/2}$. Write $Z=\{n_1<\cdots<n_h\}$, and pick an integer $m_j\in A\cap E_{n_j}$ 
for each $j=1,\ldots,h$. Then $B=\{m_1,m_3,\ldots\}=:\{\mu_1,\mu_2,\ldots\}$ is quasi-independent because 
we have $\mu_{j+1}/\mu_j\geq 2$. Moreover $\vert B\vert\geq h/2\geq(1/2)\,\vert A\vert^{1/2}$.\par
\medskip

It is now easy to see that $E$ has the required properties. Indeed, it follows from (4), (7), and from 
F. Lust-Piquard's Theorem~\ref{theo 1.9} that $E$ has a positive upper density in $\Lambda$. That it is 
CUC follows from (8) and from Proposition~\ref{prop 1.2}. That it is $\Lambda(q)$ for all $q<+\infty$ follows 
from (8) and from Proposition~\ref{prop 1.1}. The fact that $E$ is $4/3$-Rider follows from (a) in 
Proposition~\ref{prop 1.6}. Indeed, if $\eps(p)=2/p-1$, then $\eps(4/3)=1/2$.\par

Finally, let $N$ be a large integer, and $n$ such that $2^n\leq N<2^{n+1}$. Then 
\begin{align*}
\vert E\cap[1,N]\vert 
& \geq\vert E_{n_0+1}\vert+\cdots+\vert E_{n-1}\vert
\geq \frac{c\tau}{4}\big[(n_0+1)+\cdots+(n-1)\big] \\
& \geq d_0 n^2 \geq d_1(\log N)^2
\end{align*}
where $d_0,d_1$ are positive constants. If now $E$ is a $p$-Rider set, we have the mesh condition 
$\vert E\cap[1,N]\vert = O\,\big((\log N)^{p/(2-p)}\big)$. This requires $2\leq p/(2-p)$, that is 
$p\geq 4/3$. And this ends the proof of Theorem~\ref{theo 2.7}.\hfill\qed
\bigskip

\noindent{\bf Remark.} The third-named author proved the following (\cite{Luis2}, Lema 2.4) (which is actually 
implicitly already contained in \cite{Pisier48}, Lemme 7.2, Th\'eor\`eme 7.1, and Th\'eor\`eme 2.3 (iv)):\par
\begin{itemize}
\item[($\ast$)] For every finite subset $A\subseteq\Z$, there exists a quasi-independent subset $B\subseteq A$ 
such that $\vert B\vert\geq\delta (\vert A\vert / \psi_A)^2$, where $\delta>0$ is a numerical constant.\par
\end{itemize}

On the other hand, G. Pisier (\cite{Pisier2}, Lemme 5.2) proved:
\begin{equation}\label{P1}
\esp\Big\Vert \sum_k a_k r_k e_k\big\Vert_\Psi
\leq C \Big(\sum_k\vert a_k\vert^2\Big)^{1/2} \tag*{\rm (1)}
\end{equation}
where $C$ is a numerical constant, $(r_k)_k$ is the Rademacher sequence, and $\Vert\;\Vert_\Psi$ is the Orlicz 
space associated to $\Psi(x)=\e^{x^2}-1$.\par
Taking our selectors $\eps_k$ with mean $\delta_k=c\,n/2^n$ for $k\in I_n$, standard symmetrization and 
centering arguments give:
\begin{equation}\label{P2}
\esp \Big\Vert\sum_{k\in I_n}\eps_k e_k\Big\Vert_\Psi\leq C\sqrt n\,. \tag*{\rm (2)}
\end{equation}

In other terms, we have, in view of Lemma~\ref{lemme 2.9}:
\begin{equation}\label{P3}
\esp(\psi_{\Lambda_n})\leq C\sqrt n\leq C'\vert \Lambda_n\vert^{1/2}. \tag*{\rm (3)}
\end{equation}

If we could prove a concentration inequality, variant of Lemma~\ref{lemme 1.3}, then this variant and the 
Borel-Cantelli lemma would imply from (3) that:
\begin{equation}\label{P4}
\text{Almost surely } \psi_{\Lambda_n}\leq C''\vert\Lambda_n\vert^{1/2} \text{ for $n$ large enough.} 
\tag*{\rm (4)}
\end{equation}

We could then combine ($\ast$) and (4) directly to obtain the following alternative proof of 
Theorem~\ref{theo 2.7}. Select $\omega\in\Omega$ such that $\Lambda$ is strongly ergodic, with 
$\vert \Lambda_n\vert\geq c\,n/2$, and $\psi_{\Lambda_n}\leq C''\vert \Lambda_n\vert^{1/2}$; take then a 
quasi-independent set $E_n\subseteq\Lambda_n$ of size
$$\vert E_n\vert\geq\delta\left(\frac{\vert\Lambda_n\vert}{\psi_{\Lambda_n}}\right)^2 
\geq \delta C''^{-2}\vert\Lambda_n\vert\geq \delta'n\,;$$
the set $E=\bigcup\limits_n E_n$ then has the required properties.\par
\bigskip

To end this section, we consider the case $p>4/3$. We cannot keep the property of uniform convergence ({\it CUC}), 
nor that of being $q$-Sidon stated in Theorem~\ref{theo 2.5}. We do not know whether this is only due to the method. 
But being $p$-Rider with $p>4/3$ might be a rather weak condition (see \cite{LLQR} and 
\cite{Lefevre-Luis}).\par
\medskip

\begin{theoreme}\label{theo 2.10}
For every $4/3\leq p<2$ there exists a set $\Lambda$ of integers which is $p$-Rider, but is not $q$-Rider for 
$q<p$ and which is $\Lambda(q)$ for every $q<+\infty$, but which is uniformly distributed (so in particular 
dense in the Bohr group, and ${\cal C}_\Lambda$ contains $c_0$).
\end{theoreme}

The proof is essentially the same as in Theorem~\ref{theo 2.2}, except that we take selectors $\eps_k$ of mean
$$ \delta_k=c\, \frac{(\log k)^\alpha}{k(\log\log k)^{\alpha+1}} \quad \text{for }  k\geq 1,$$
where $\alpha=2(p-1) / (2-p) \geq 1$, and replace $M_n=n^n$ by the smallest integer $\geq n^{\beta n}$, with 
$\beta$ any number $>\alpha$ (for instance $\beta=\alpha+1$), which we call again $M_n$. The estimate:
$$\prob\big(\Omega_n(M_n)\big)\leq 2\, \frac{c^n}{n^n}\, \frac{(\log M_n)^{n(\alpha+1)}}{M_n}$$
still holds, and now gives:
$$\prob\big(\Omega_n(M_n)\big)\leq C'\,^n\,\frac{(\log n)^{n(\alpha+1)}} {n^{(\beta-\alpha)n}}\,\cdot$$
Then easy computations show that:\par
($\ast$) Almost surely $\vert \Lambda_{M_n}\vert\sim(n\log n)^{\alpha+1}$ for $n$ sufficiently large;\par
($\ast\ast$) Almost surely $\vert\Lambda'_n\vert\sim n^\alpha(\log n)^{\alpha+1}$ for $n$ sufficiently 
large.\par

\noindent
Property ($\ast$) guaranties that $\Lambda(\omega)$ will still be almost surely $p$-Rider, and ($\ast\ast$) with 
the mesh condition implies that $\Lambda$ is not $q$-Rider for $q<p$.\par
The $\Lambda(q)$ property cannot be obtained by the Littlewood-Paley method, but follows from 
\cite{Neuwirth}, Theorem 4.7.\hfill\qed\par


\section{Large thin sets in prescribed sets of integers}

In this section, we start from a prescribed set 
$\Lambda_0 =\{\lambda_1<\lambda_2<\ldots<\lambda_N<\ldots\}$ of positive integers, and randomly 
construct a thin set $\Lambda$ inside $\Lambda_0$ in the following way. We still have our selectors 
$\eps_1,\ldots,\eps_N,\ldots$ of respective means $\delta_1,\ldots,\delta_N,\ldots$. This time, however, 
we set
$$\Lambda=\Lambda(\omega)=\{\lambda_j\in\Lambda_0\,;\ \eps_j(\omega)=1\},$$
{\it i.e.} we select randomly some of the $\lambda_j$'s, and ignore the other integers. Such constructions have been 
made previously by S. Neuwirth (\cite{Neuwirth}).\par \smallskip

We always assume that $\Lambda_0$ is ergodic, namely that
$$A_{\Lambda_0,N}(t)=N^{-1}\big(\e^{i\lambda_1t}+\cdots+\e^{i\lambda_Nt}\big)
\mathop{\longrightarrow}_{N\to+\infty} l(t)\,,\hskip 5mm \forall t\in\T.$$\par\medskip

In this context, we have the following theorem, which extends Bourgain's Theorem~\ref{theo 1.10}.

\begin{theoreme}[\cite{Neuwirth}, Th. 5.4]\label{theo 3.1}
Let $\Lambda_0$ be an ergodic (resp. strongly ergodic, resp. uniformly distributed) set of positive integers, and let 
$\eps_1,\ldots,\eps_N,\ldots$ be selectors with respective expectation $\delta_1,\ldots,\delta_N,\ldots$ with 
$(\delta_n)_{n\geq1}$ decreasing. Assume that 
$\sigma_N/\log\lambda_N\mathop{\longrightarrow}\limits_{N\to+\infty} +\infty$, 
where $\sigma_N=\delta_1+\cdots+\delta_N$. Then, almost surely, the set $\Lambda$ is ergodic (resp. strongly 
ergodic, resp. uniformly distributed). More precisely, if 
$A_{\Lambda_0,N}(t)\mathop{\longrightarrow}\limits_{N\to+\infty} l(t)$,
we have, almost surely, with $\Lambda_N=\Lambda\cap\{\lambda_1,\ldots,\lambda_N\}$,
$$A_N(t)= \frac{1}{\vert \Lambda_N\vert}\sum_{n\in\Lambda_N}e_n(t)
\mathop{\longrightarrow}_{N\to+\infty}l(t)\,,\hskip 5mm\forall t\in\T.$$
\end{theoreme}

We sketch the proof. First, we require

\begin{lemme}\label{lemme 3.2}
Let $\eps_1,\ldots,\eps_N$ be selectors of respective expectations $\delta_1,\ldots,\delta_N$. Setting 
$\sigma_N=\delta_1+\cdots+\delta_N$, one has the following inequality:
$$\prob\big(\big\Vert{\scriptstyle\sum\limits_{k=1}^N}  
(\eps_k-\delta_k)e_{\lambda_k}\big\Vert_\infty > 15\sqrt{\sigma_N\log \lambda_N}\big) 
\leq 8 / N^2\,,$$
provided that $\sigma_N\geq 25\log \lambda_N\,.$
\end{lemme}

\noindent{\bf Proof.} Set $Q=\sum\limits_{k=1}^N (\eps_k-\delta_k)e_{\lambda_k}$. For fixed $t\in\R$, one 
has $Q(t)=\sum\limits_{k=1}^N X_k$, where $X_k=e_{\lambda_k}(t)(\eps_k-\delta_k)$. The $X_k$'s are 
independent, bounded by $1$, and centered complex random variables; so, letting 
$t_N=5\sqrt{\sigma_N\log \lambda_N}$, and using Lemma~\ref{lemme 1.3}, we get
\begin{align*}
\prob(\Vert Q \Vert_\infty>3t_N) 
&\leq \prob(\sup_{t\in F_N}\vert Q(t)\vert > t_N)
\leq\sum_{t\in F_N}\prob(\vert Q(t)\vert>t_N)\cr
&\leq 4\vert F_N\vert\,\exp\big(- \frac{t_N^2}{8\sigma_N}\big) =
8\lambda_N^{1-25/8}\leq N^{1-25/8}\leq 8N^{-2},
\end{align*}
where $F_N=\big\{ j\pi / \lambda_N\,;\ 0\leq j\leq 2\lambda_N-1\,\big\}$ is the set of the 
$(2\lambda_N)^{th}$ roots of unity, and where the first inequality follows from Bernstein 
inequality (see \cite{Frappier}).\hfill\qed
\bigskip

\noindent{\bf Proof of Theorem~\ref{theo 3.1}.} 
Notice first that
$$\frac{1}{\sigma_N}\sum_{n=1}^N \delta_n e_{\lambda_n}(t)
\mathop{\longrightarrow}_{N\to+\infty}l(t)\,,\hskip 5mm \forall t\in\T\,.$$
\par

In fact, set $E_n=e_{\lambda_1}(t)+\cdots+e_{\lambda_n}(t)$ and $l=l(t)$. Since $(\delta_n)_n$ is 
nonincreasing,  two Abel's partial summations give:
\begin{align*}
\sum_{n=1}^N \delta_ne_{\lambda_n}(t) 
& =\sum_{n=1}^{N-1}(\delta_n-\delta_{n+1})E_n+\delta_N E_N \cr
&=\sum_{n=1}^{N-1} nl (\delta_n-\delta_{n+1}) + Nl\delta_N
+ o\,\Big(\sum_{n=1}^{N-1} n (\delta_n-\delta_{n+1})+N\delta_N\Big) \cr
& = l \sigma_N+ o\, (\sigma_N).
\end{align*}

Setting $Q_N=\sum\limits_{n=1}^N(\eps_n-\delta_n)e_{\lambda_n}$, we have:
$$\Big\Vert A_N- \frac{1}{\sigma_N}
\sum_{n=1}^N\delta_n e_{\lambda_n}\Big\Vert_\infty
\leq \frac{2\Vert Q_N\Vert_\infty}{\sigma_N}\, \raise 1pt \hbox{,}$$
since
\begin{align*}
\Big\Vert \frac{1}{ S_N}\sum_{n=1}^N \eps_n e_{\lambda_n} 
& - \frac{1}{\sigma_N}
\sum_{n=1}^N \delta_n e_{\lambda_n}\Big\Vert_\infty\cr
&\leq \left\vert \frac{1}{S_N} - \frac{1}{\sigma_N}\right\vert\,
\Big\Vert\sum_{n=1}^N \eps_n e_{\lambda_n}\Big\Vert_\infty
+ \frac{1}{\sigma_N} \Big\Vert\sum_{n=1}^N (\eps_n-\delta_n) e_{\lambda_n}\Big\Vert_\infty\cr
&\leq \frac{\vert S_N-\sigma_N\vert}{\sigma_N} + \frac{\Vert Q_N\Vert_\infty}{\sigma_N} =
\frac{\vert Q_N(0)\vert+\Vert Q_N\Vert_\infty }{\sigma_N}
\leq \frac{2\Vert Q_N\Vert_\infty}{\sigma_N}\,\cdot
\end{align*}

Now, Lemma~\ref{lemme 3.2} gives:
$$\prob\big(\Vert Q_N\Vert_\infty>15\sqrt{\sigma_N\log\lambda_N}\big)\leq
8\lambda_N^{-2}\leq 8N^{-2}$$
if $\sigma_N\geq 25\log\lambda_N$; so we get, by the Borel-Cantelli lemma, 
$$\frac{\Vert Q_N\Vert_\infty}{\sigma_N} =
O\,\Big(\sqrt{ \frac{\log\lambda_N}{\sigma_N}}\Big)$$
almost surely. In view of the hypothesis, we have:
$$\big\Vert A_N- \frac{1}{\sigma_N} \sum_{n=1}^N\delta_n e_{\lambda_n}\big\Vert_\infty
\mathop{\longrightarrow}_{N\to+\infty}0\hskip 5mm \hbox{almost surely;}$$
and so, almost surely $A_N(t)\mathop{\longrightarrow}\limits_{N\to+\infty} l(t)$ for each $t$, which is the 
desired conclusion.\hfill\qed
\bigskip

{\large \bf 3.2 Regularity}
\medskip

Let $I$ be a finite interval of $\N^\ast$ and $\nu(I)=\vert \Lambda_0\cap I\vert$ be the number of indices $n$ 
for which $\lambda_n\in I$. In the sequel, we assume that $\Lambda_0$ has the following regularity property:
\par
\bigskip

{\it There \hfill exists \hfill a \hfill continuous \hfill eventually \hfill strictly \hfill increasing \hfill function\\ 
$\phi \colon \,]\,0,+\infty[\,\to\, ]\,0,+\infty[$ such that}:
\begin{equation}\label{eq 3.1}
\frac{\nu([N, 2N[)}{\phi(N)} \mathop{\longrightarrow}_{N\to+\infty} 1
\end{equation}
{\it and}:
\begin{equation}\label{eq 3.2}
\frac{\phi(2x)}{\phi(x)}\mathop{\longrightarrow}_{x\to+\infty} l > 1.
\end{equation}
Note that $l\leq2$, since $\nu([2^k,2^{k+1}[\leq 2^k$ implies that 
$(1-\eps)^{k-k_0}l^{k-k_0}\phi(2^{k-k_0})\leq \phi(2^k)\leq (1+\eps)\,2^k$.
\par
We say that $\Lambda_0$ is {\it regular} if these properties hold.\par\smallskip

They are obviously verified when $\lambda_n=n^s$, and also, by the Prime Number Theorem, when 
$\lambda_n=p_n$, with $\phi(x)=x/\log x$.\par

It is easy to see that \eqref{eq 3.1} and \eqref{eq 3.2} imply that $\Lambda_0$ has a polynomial growth, namely 
that there exist two constants, $a,d>0$ such that:
\begin{equation}\label{eq 3.3}
\nu([1,k]) \geq a\,k^d.
\end{equation}
(or, equivalently, $\lambda_N\leq a'N^{1/d}$).\par 
It follows that the condition 
$\sigma_N/\log \lambda_N \mathop{\longrightarrow}\limits_{N\to+\infty}+\infty$ of Theorem~\ref{theo 3.1}  
reduces then to the previous condition 
$\sigma_N/\log N \mathop{\longrightarrow}\limits_{N\to+\infty}+\infty$ of Theorem~\ref{theo 1.10}.\par
Moreover $\Lambda_0$ satisfies:
\begin{equation}\label{eq 3.4}
\lambda_{8n}\geq 2\lambda_n\hskip 10mm\hbox{\it for $n\geq 1$ large enough.}
\end{equation}
Indeed, if $\nu([1,2^{k-1}[)<n\leq\nu([1,2^k[)$, then $2\lambda_n\leq 2^{k+1}$ and it suffices to show that 
$\nu([2^{k-1},2^{k+1}[)\leq 7n$. But 
\begin{align*}
\nu([2^{k-1},2^{k+1}[)&\leq (1+\eps)\big(\phi(2^{k-1})+\phi(2^k)\big)
\leq (1+\eps)^2(l^2+l)\phi(2^{k-2})\cr
&\leq (1+\eps)^3(l^2+l)\nu([2^{k-2},2^{k-1}[)
\leq (1+\eps)^3(l^2+l)\nu([1,2^{k-1}[)\cr
&\leq 7\nu([1,2^{k-1}[)\leq 7n
\end{align*}
for $\eps>0$ small enough, and $n$ large enough.\par\medskip

As in Section~\ref{section 2}, we restrict ourselves to selectors with mean $\delta_n=\alpha_n/n$, where 
$(\delta_n)_n$ decreases to $0$, and $(\alpha_n)_n$ tends to infinity, and moreover, except in the last theorem, 
$(\alpha_n)_n$ increases.\par

The following lemma is quite similar to Lemma~\ref{lemme 2.1}. We indicate some changes which are needed, and 
how the regularity occurs.\par

\begin{lemme}\label{lemme 3.3}
Let $s\geq2$ and $M$ be integers and let
$$\Omega_s(M)=\{\omega\in\Omega\,;\ \Lambda(\omega)\cap[\lambda_M, +\infty[
\hbox{ contains at least a relation of length $s$ }\}.$$\par
We have, for $s$ large enough,
$$\prob\big(\Omega_s(M)\big)\leq \frac{(16\e)^s}{s^s} \sum_{j>M}\delta_j^2\sigma_j^{s-2}.$$
\end{lemme}

\noindent{\bf Proof.} As in the proof of Lemma~\ref{lemme 2.1}, we write 
$\Omega_s(M)=\bigcup\limits_{l\geq M+s-1}\Delta_l$, where $\Delta_l$ is defined by
\begin{align*}
\Delta_l =\{\omega\,;\ \Lambda(\omega)\cap[\lambda_M, +\infty[ 
\text { contains at least } 
&\text{a relation of length $s$} \\ 
& \text{and with greatest term $\lambda_l $}\}.
\end{align*}

It suffices to show that 
$$\prob(\Delta_l)\leq \frac{8^s 2^{s-2}}{(s-2)!}\,\delta_l^2\sigma_l^{s-2}.$$
The proof proceeds as in Lemma~\ref{lemme 2.1}, replacing $i_1,\ldots,i_{s-1}$ and $l$ by 
$\lambda_{i_1},\ldots,\lambda_{i_{s-1}}$ and $\lambda_l$ respectively. The relation ($\ast\ast$) gives 
$\lambda_{i_{s-1}}\geq \lambda_l/s$. The regularity appears now to say that $i_{s-1}\geq l/8^s$. Indeed, 
otherwise, by \eqref{eq 3.4}, we should have, for $s$ large enough,
$$\lambda_l>\lambda_{8^s i_{s-1}}\geq 2^s\lambda_{i_{s-1}} \geq 2^s \frac{\lambda_l}{s} >\lambda_l.$$\par 

This gives the lemma since $(\alpha_n)_n$ increases:
\begin{equation}
\delta_{i_{s-1}} = \frac{\alpha_{i_{s-1}}}{i_{s-1}} \leq \frac{\alpha_l}{i_{s-1}}=
\frac{\alpha_l}{ l} \frac{l}{i_{s-1}} \leq 8^s\delta_l\,. \tag*{$\square$}
\end{equation}

Since this basic lemma still holds for random subsets of prescribed sets $\Lambda_0$, the first main theorems of 
Section~\ref{section 2} still hold and their proofs requires only minor modifications because of 
Theorem~\ref{theo 3.1}. We therefore content ourselves with stating them.

\begin{theoreme}
Let $\Lambda_0$ be a regular, strongly ergodic set of positive integers. There exists a set 
$\Lambda\subseteq\Lambda_0$ which is:
\begin{itemize}
\item[{\rm (1)}] $p$-Sidon for all $p>1$, $\Lambda(q)$ for all $q<+\infty$, {\rm CUC}, but which is:\par
\item[{\rm (2)}] strongly ergodic (in particular, ${\cal C}_\Lambda$ contains $c_0$ and $\Lambda$ is not a 
Rosenthal set).
\end{itemize}
\end{theoreme}

\begin{theoreme}
Let $\Lambda_0$ be as in the previous theorem, and let $1<p<4/3$. Then, there exists a set 
$\Lambda\subseteq\Lambda_0$ which is:
\begin{itemize}
\item[{\rm (1)}] strongly ergodic (in particular, ${\cal C}_\Lambda$ contains $c_0$ and so $\Lambda$ is 
not a Rosenthal set), but which is:\par
\item[{\rm (2)}] a {\rm CUC}-set, $\Lambda(q)$ for all $q<+\infty$, and\par
\qquad {\rm (a)} is $p$-Rider, but is not $q$-Rider for $q<p$,\par
\qquad {\rm (b)} is $q$-Sidon for all $q>p/(2-p)$.
\end{itemize}
\end{theoreme}

\begin{theoreme}
Same as in the previous theorem, but instead of property {\rm (a)}:\par
{\rm (a')} $\Lambda$ is $q$-Rider for every $q>p$, but is not  $p$-Rider.
\end{theoreme}

\begin{theoreme}
Let $\Lambda_0=\{\lambda_1,\ldots\}$ be a regular, strongly ergodic set. Then, there exists a set 
$\Lambda\subseteq\Lambda_0$ which is strongly ergodic and contains a set $E$ which\par
\begin{itemize}
\item[{\rm (1)}] has a positive upper density in $\Lambda(\omega)$ (so in particular, ${\cal C}_E$ contains 
$c_0$ and $E$ is not a Rosenthal set), and 
\par
\item[{\rm (2)}] is a {\rm CUC}-set, is $4/3$-Rider, but not $q$-Rider for $q<4/3$, and is a $\Lambda(q)$-set 
for all $q<+\infty$; more precisely, for all $q>2$, we have $\Vert f\Vert_q\leq Cq^2\Vert f\Vert_2$ for all 
$f\in{\cal P}_E$, where $C>0$ is a numerical constant.
\end{itemize}
\end{theoreme}

The proof is the same as that of Theorem~\ref{theo 2.7}, so we omit it. We merily note the following facts.\par

The sequence $(\delta_k)_k$ is eventually decreasing. Indeed, for $n\geq n_\eps$, we have, by the regularity 
conditions \eqref{eq 3.1} and \eqref{eq 3.2}, if $\eps>0$ is chosen so that $(1-\eps)^2 l\geq 1$,
$$\nu_{n+1}\geq (1-\eps)\phi(2^{n+1})\geq (1-\eps)^2 l\,\phi(2^n) \geq(1-\eps)^3 l\,\nu_n\geq\nu_n\,.$$\par

Next, \eqref{eq 3.1} implies that, for some constant $\alpha>0$, and for $2^q<N\leq 2^{q+1}$,
\begin{align*}
\sigma_N =\delta_1+\cdots+\delta_N 
& \geq \sum_{n=1}^q\Big(\sum_{k\in I_n}\delta_k\Big) = c\,\sum_{n=1}^q \log \nu_n \\
& \geq c\alpha\,\sum_{n=1}^q n\geq c(\alpha/2)q^2.
\end{align*}
Since $\Lambda_0$ has polynomial growth$\,:$ $\lambda_N= O\,(N^{1/d})$, we have 
$\log \lambda_N\leq \lambda_{2^{q+1}}= O\,(q)$. It follows that:
$$\sigma_N/\log \lambda_N
\mathop{\longrightarrow}\limits_{N\to+\infty}+\infty\,.$$\par

Finally, we have to replace the parameter $\psi_A$ in the proof of Theorem~\ref{theo 2.7} by:
$$\psi'_A=\sup_{p\geq 2} \frac{\Vert e'_A\Vert_p}{\sqrt p},\hskip 5mm
\hbox{where } e'_A=\sum_{\lambda_k\in A}e_{\lambda_k}\,.$$
Since we have, for any interval $I\,:$
$\psi'_I\leq C\,\nu(I) / \sqrt{\log \nu(I)}$, the rest of the proof will then work with no essential change.\hfill\qed
\par\bigskip

\noindent{\bf Remark.} Consider the $\psi$-parameter associated to the squares, that is:
$$\psi'_N=\sup_{q\geq 2} \frac{\Vert S_N\Vert_q}{\sqrt q}\,  \raise 1pt\hbox{,}$$
where $S_N(x)=\sum\limits_{n=1}^N \e^{in^2x}$. It follows from results of 
Zalcwasser (\cite{Zalc}), that we have very precise estimates on $\Vert S_N\Vert_q$: there exist numerical 
constants $C_1,C_2>0$ such that:
$$C_1N^{1-2/q}\leq \Vert S_N\Vert_q\leq C_2N^{1-2/q}$$
whenever $q\geq 5$ and $N\geq 1$ (when $q$ is near 4, a logarithmic factor $(\log N)^{1/q}$ should be added in 
the upper estimate). Therefore, the {\it a priori} crude estimate used in the proof of Theorem~\ref{theo 2.7} is, at 
least for the squares, optimal, as it is for the set of all the positive integers.\par



\begin{thebibliography}{99}

\bibitem{Arki-Oskol} G.I. Arkipov and K.I. Oskolkov, 
On a special trigonometric series and its applications, 
Math. Sbornik 134 (1987), 145--155.

\bibitem{Blei1} R. Blei, 
Sidon partitions and $p$-Sidon sets, 
Pacific J. Math. 65 (1976), 307--313.

\bibitem{Blei2} R. Blei, 
Fractional cartesian products of sets, 
Ann. Inst. Fourier 29 (1979), 79--105.

\bibitem{Blei3} R. Blei, 
Combinatorial dimension and certain norms in Harmonic Analysis, 
Amer. J. Math. 106 (1984), 847--887.

\bibitem{Blei-K} R. Blei and T.W. K\"orner, 
Combinatorial dimension and random sets, 
Israel J. Math. 47 (1984), 65--74.

\bibitem{Blum} J.R. Blum, B. Eisenberg and L.-S. Hahn, 
Ergodic theory and the measure of sets in the Bohr group, 
Acta Scient. Math. Szeged 34 (1973), 17--24.

\bibitem{Bourgain7} J. Bourgain, 
Propri\'et\'es de d\'ecomposition pour les ensembles de Sidon, 
Bull. Soc. Math. France 111 (1983), 421-428.

\bibitem{Bourgain8} J. Bourgain, 
Sidon sets and Riesz products, 
Ann. Inst. Fourier 35 (1985), 137--148.

\bibitem{Bourgain9} J. Bourgain, 
Subspaces of $L^\infty_N$, arithmetical diameter and Sidon sets, 
Probability in Banach Spaces V, Proceed. Medford 1984, 
Lecture Notes in Math. 1153 (1985), 96--127.

\bibitem{Bourgain10} J. Bourgain, 
On the maximal ergodic theorem for certain subsets of the integers, 
Israel J. Math. 61 (1988), 39--72.

\bibitem{Bourgain11} J. Bourgain, 
On the behavior of the constant in the Littlewood-Paley inequality, 
Geometric aspects of functional analysis,  
Isr. Semin., GAFA, Isr. 1987-88, Lecture Notes in Math. 1376 (1989), 202-208.

\bibitem{Bourgain12} J. Bourgain, 
Bounded orthogonal sets and the $\Lambda(p)$-set problem, 
Acta Math. 162 (1989), 227--246.

\bibitem{Bourgain-Milman} J. Bourgain and V. Milman, 
Dichotomie du cotype pour les espaces invariants, 
C.R.A.S. Paris 300 (1985), 263--266.

\bibitem{Bo-Pyt} M. Bo$\dot {\rm z}$ejko and T. Pytlik, 
Some types of lacunary Fourier series, 
Colloq. Math. 25 (1972), 117--124.

\bibitem{Myriam} M. D\'echamps-Gondim, 
Sur les compacts associ\'es aux ensembles lacunaires, les ensembles de Sidon et quelques probl\`emes ouverts,  
Publications Math\'ematiques d'Orsay 84-01 (1984), expos\'e 7.

\bibitem{Diestel} J. Diestel, 
Sequences and Series in Banach Spaces, 
Graduate Texts in Math. 92, Springer-Verlag (1984).

\bibitem{Drury} S. Drury, 
Sur les ensembles de Sidon, 
C.R.A.S. Paris 271 (1970), 162--163.

\bibitem{Ed-Ross} R.E. Edwards and K.A. Ross, 
$p$-Sidon sets, 
J. Funct. Anal. 15 (1974), 404--427.

\bibitem{Figa} A. Fig\`a--Talamanca, 
An example in the theory of lacunary Fourier series, 
Boll. Unione Matem. Ital. 3 (1970), 375--378.

\bibitem{Fournier} J. Fournier, 
Two UC-sets whose union is not a UC-set, 
Proc. Amer. Math. Soc. 84 (1982), 69--72.

\bibitem{Fournier-P} J. Fournier and L. Pigno, 
Analytic and arithmetic properties of thin sets, 
Pacific J. Math. 105 (1983), 115--141.

\bibitem{Frappier} C. Frappier, Q.I. Rahman and St. Ruscheweyh, 
New inequalities for polynomials, 
Trans. Amer. Math. Soc. 288 (1985), 69-99.

\bibitem{Johnson} G.W. Johnson, 
Theorems on lacunary sets, especially $p$-Sidon sets, 
Studia Math. 58 (1976), 209--221.

\bibitem{Johnson-W} G.W. Johnson and G.S. Woodward, 
On $p$-Sidon sets, 
Indiana Univ. Math. J. 24 (1974), 161--167.

\bibitem{Kahane} J.--P. Kahane, 
Some random series of functions, Second ed., 
Cambridge Univ. Press, Cambridge (1985).

\bibitem{Kashin-T} B. Kashin and L. Tzafriri, 
On random sets of uniform convergence, 
Math. Notes 54 (1993), 677--687.

\bibitem{Katz} Y. Katznelson, 
Suites al\'eatoires d'entiers, 
Lecture Notes in Math. 336, Springer-Verlag Berlin (1973), 148-152.

\bibitem{Katz-Mall} Y. Katznelson and  P. Malliavin, 
V\'erification statistique de la conjecture de la dichotomie sur une classe d'alg\`ebres de restriction,  
C.R.A.S. Paris 262 (1966), 490--492.

\bibitem{Konyagin} S. Konyagin, 
On divergence of trigonometric Fourier series everywhere, 
C.R.A.S. Paris 329 (1999), 693--697.

\bibitem{Korner} T.W. K\"orner, 
Fourier analysis, 
Cambridge University Press (1988).

\bibitem{Kwap-P} S. Kwapien and A. Pe{\l}czy\'nski, 
Absolutely summing operators and translation invariant spaces of functions on compact abelian groups,  
Math. Nachrichten 94 (1980), 303--340.

\bibitem{Ledoux-Ta} M. Ledoux and M. Talagrand, 
Probability in Banach spaces, 
Ergeb. Math. 23 Springer-Verlag (1991).

\bibitem{Lefevre1} P. Lef\`evre, 
Sur les ensembles de convergence uniforme, 
Publ. Math. d'Orsay 94--24 (1994), 1--70.

\bibitem{Lefevre2} P. Lef\`evre, 
Measures and lacunarity sets, 
Studia Math. 133 (1999), 145--161.

\bibitem{LLQR} P. Lef\`evre, D. Li, H. Queff\'elec, and L. Rodr\'{\i}guez-Piazza,  
Lacunary sets and function spaces with finite cotype, ({\it submitted})

\bibitem{Lefevre-Luis} P. Lef\`evre and L. Rodr\'{\i}guez-Piazza, 
$p$-Rider sets are $q$-Sidon sets, ({\it submitted})

\bibitem{Li} D. Li, 
A remark about $\Lambda(p)$-sets and Rosenthal sets, 
Proc. Amer. Math. Soc. 126 (1998) 3329--3333.

\bibitem{L-T} J. Lindenstrauss and L. Tzafriri, 
Classical Banach Spaces I and II, 
Classics in Math., Springer (1997).

\bibitem{Lopez-Ross} J.M. Lopez and K.A. Ross, 
Sidon Sets, 
Marcel Dekker 13 (1975).

\bibitem{Lust1} F. Lust, 
Produits tensoriels injectifs d'espaces de Sidon, 
Colloq. Math. 32 (1975), 285--289.

\bibitem{Lust41} F. Lust--Piquard, 
Propri\'et\'es g\'eom\'etriques des sous-espaces invariants par translation de $L^1(G)$ et ${\cal C}(G)$,  
S\'emin. G\'eom. Espaces Banach, Ecole Polytechnique, Paris (1977-78), Expos\'e n~$^\circ$26.

\bibitem{Lust42} F. Lust--Piquard, 
Bohr local properties of ${\cal C}_\Lambda(\T)$, 
Colloq. Math. 58 (1989), 29--38.

\bibitem{Neuwirth} S. Neuwirth, 
Random constructions inside lacunary sets,  
Annales Inst. Fourier 49 (1999), 1853--1867.

\bibitem{Oskol} K.I. Oskolkov, 
On spectra of uniform convergence, 
Soviet Math. Dokl. 33 (1986), 616--620.

\bibitem{Pedem} L. Pedemonte, 
Sets of uniform convergence, 
Colloq. Math. 33 (1975), 123--132.

\bibitem{Pisier1} G Pisier, 
Ensembles de Sidon et espaces de cotype 2,  
S\'eminaire sur la g\'eom\'etrie des espaces de Banach 1977--1978, Ecole Polytechnique, Paris (1978), expos\'e 14.

\bibitem{Pisier2} G. Pisier, 
Sur l'espace de Banach des s\'eries de Fourier al\'eatoires presque s\^urement continues, 
S\'eminaire sur la g\'eom\'etrie des espaces de Banach 1977--1978, Ecole Polytechnique, 
Paris (1978), expos\'es 17--18.

\bibitem{Pisier48} G. Pisier, 
De nouvelles caract\'erisations des ensembles de Sidon, 
Math. Anal. and Applic., Part B, Advances in Math. Suppl. Studies, Vol 7B (1981), 685--726.

\bibitem{Rider1} D. Rider, 
Gap series on groups and spheres, 
Canad. J. Math. 18 (1966), 389--398.

\bibitem{Rider2} D. Rider, 
Randomly continuous functions and Sidon sets, 
Duke Math. J. 42 (1975), 759--764.

\bibitem{Luis1} L. Rodr\'\i guez--Piazza, 
Caract\'erisation des ensembles $p$-Sidon p.s., 
C.R.A.S. Paris 305 (1987), 237--240.

\bibitem{Luis2} L. Rodr\'\i guez--Piazza, 
Rango y propiedades de medidas vectoriales. Conjuntos $p$-Sidon p.s., 
Thesis, Universidad de Sevilla (1991).

\bibitem{Rosenthal} H.P. Rosenthal, 
On Trigonometric Series Associated with Weak$^\ast$ Closed Subspaces of Continuous Functions, 
Journ. Math. Mech. 17 (1967), 485-490.

\bibitem{Rudin} W. Rudin, 
Trigonometric Series with Gaps, 
Journal of Math. and Mech. 9 (1960), 203--227.

\bibitem{Singer} I. Singer, 
Bases in Banach Spaces I, 
Springer Verlag (1970).

\bibitem{Soardi-T} P.M. Soardi and G. Travaglini, 
On sets of completely uniform convergence, 
Colloq. Math. 45 (1981), 317--320.

\bibitem{Trava} G. Travaglini, 
Some properties of UC-sets, 
Boll. Unione Matem. Ital. 15 (1978), 272--284.

\bibitem{Varo} N.T. Varopoulos, 
Une remarque sur les ensembles de Helson,  
Duke Math. J. 43 (1976), 387--390.

\bibitem{Wojt} P. Wojtaszczyk, 
Banach Spaces for Analysts, 
Cambridge University Press (1991).

\bibitem{Woodward} G. S. Woodward, 
$p$-Sidon Sets and a Uniform Property, 
Indiana Univ. Math. Journal 25 (1976), 9951--1003.

\bibitem{Zalc} Z. Zalcwasser, 
Polyn\^omes associ\'es aux fonctions modulaires $\vartheta$, 
Studia Math. 7 (1938), 16--35.

\bibitem{Zyg} A. Zygmund, 
Trigonometric Series, Second Ed., Vol. I \& II, 
Cambridge Math. Library, Cambridge Univ. Press (1993). 

\end{thebibliography}
\end{document}